\documentclass{amsart}

\usepackage{amsmath} 
\usepackage{amssymb}
\usepackage{mathrsfs}

\newtheorem{theorem}{Theorem}[section] 
\newtheorem{claim}[theorem]{Claim}

\newtheorem{observation}[theorem]{Observation}

\theoremstyle{definition}
\newtheorem{definition}[theorem]{Definition}
\newtheorem{convention}[theorem]{Convention}

\newtheorem{discussion}[theorem]{Discussion}

\theoremstyle{remark}
\newtheorem{remark}[theorem]{Remark}

\newtheorem{fact}[theorem]{Fact}
\newtheorem{notation}[theorem]{Notation}

\newcommand{\pcf}{{\rm pcf}}

\newcommand{\dual}{{\rm dual}}
\newcommand{\comp}{{\rm comp}}

\newcommand{\refl}{{\rm refl}}

\newcommand{\dg}{{\rm dg}}
\newcommand{\ch}{{\rm ch}}
\newcommand{\cd}{{\rm cd}}

\newcommand{\pp}{{\rm pp}}

\newcommand{\ac}{{\rm ac}}
\newcommand{\bd}{{\rm bd}}
\newcommand{\ged}{{\rm gd}}
\newcommand{\Ord}{{\rm Ord}}
\newcommand{\acc}{{\rm acc}}
\newcommand{\sch}{{\rm sch}}
\newcommand{\stat}{{\rm stat}}
\newcommand{\nacc}{{\rm nacc}}

\newcommand{\Reg}{{\rm Reg}}
\newcommand{\GCH}{{\rm GCH}}

\newcommand{\tcf}{{\rm tcf}}

\newcommand{\otp}{{\rm otp}}

\newcommand{\cg}{{\rm cg}}

\newcommand{\good}{{\rm good}}

\newcommand{\Rang}{{\rm Rang}}

\newcommand{\rest}{{\restriction}}

\newcommand{\wilog}{{\rm without loss of generality}}

\newcommand{\then}{{\underline{then}}}
\newcommand{\when}{{\underline{when}}}
\newcommand{\Then}{{\underline{Then}}}

\newcommand{\mn}{{\medskip\noindent}}
\newcommand{\sn}{{\smallskip\noindent}}

\newcommand{\cD}{{\mathscr D}}

\newcommand{\cH}{{\mathscr H}}

\newcommand{\cF}{{\mathscr F}}

\newcommand{\cP}{{\mathscr P}}

\newcommand{\cf}{{\rm cf}}

\newcount\skewfactor
\def\mathunderaccent#1#2 {\let\theaccent#1\skewfactor#2
\mathpalette\putaccentunder}
\def\putaccentunder#1#2{\oalign{$#1#2$\crcr\hidewidth
\vbox to.2ex{\hbox{$#1\skew\skewfactor\theaccent{}$}\vss}\hidewidth}}

\newenvironment{PROOF}[2][\proofname.]
{\begin{proof}[#1]\renewcommand{\qedsymbol}{\textsquare$_{\rm #2}$}}
   {\end{proof}}

\begin{document}

\title {non-reflection of the bad set for
$\check I_\theta[\lambda]$ and $\pcf$}
\author {Saharon Shelah}
\address{Einstein Institute of Mathematics\\
Edmond J. Safra Campus, Givat Ram\\
The Hebrew University of Jerusalem\\
Jerusalem, 91904, Israel\\
 and \\
 Department of Mathematics\\
 Hill Center - Busch Campus \\ 
 Rutgers, The State University of New Jersey \\
 110 Frelinghuysen Road \\
 Piscataway, NJ 08854-8019 USA}
\email{shelah@math.huji.ac.il}
\urladdr{http://shelah.logic.at}
\thanks{We would like to thank the Israel Science Foundation for
  partial support of this research. Publication 1008.\\
The author thanks Alice Leonhardt for the beautiful typing.}



\date{June 12, 2012}

\begin{abstract}
We reconsider here the 
following related pcf questions and make some advances:
\bigskip

\noindent
(Q1) \quad concerning the ideal $\check I_\kappa[\lambda]$
how much reflection do we have for the bad set
$S^{\bd}_{\lambda,\kappa} \subseteq \{\delta < \lambda:\cf(\delta) =
\kappa\}$ assuming it is well defined, (for transparency only)?
\smallskip

\noindent
(Q2) \quad are there somewhat free black boxes?

\noindent
The advances in (Q2) will be used in subsequent for constructions of
Abelian groups and modules.
\end{abstract}

\maketitle
\numberwithin{equation}{section}
\setcounter{section}{-1}
\newpage

\centerline {Anotated Content}
\bigskip

\noindent
\S0 \quad Introduction, \pageref{Introduction}
\smallskip

\S(0A) \quad Background, \pageref{Background}
\smallskip

\S(0B) \quad Results, \pageref{Results}

\hskip30pt (label y-)
\smallskip

\S(0C) \quad Quoting Definitions, \pageref{Quoting}

\hskip30pt (labels g-,m-m-m3,m6,m9,m13,m15,m52)
\bigskip

\noindent
\S1 \quad On Systems, \pageref{On}

\hskip30pt (labels b-b4, proof of \ref{y12})
\smallskip

\S(1A) \quad Existence of large members of $\check I_\theta[\lambda]$,
\pageref{Existence}
\smallskip

\S(1B) \quad Quite free witnesses of pcf-cases exists, \pageref{Quite}
\newpage

\section {Introduction} \label{Introduction}
\bigskip

\subsection {Background} \label{Background} \
\bigskip

On $\check I_\theta[\lambda]$ for $\lambda > \theta$ regular see
(Definition \ref{m13}(3) and) \cite{Sh:108}, 
\cite{Sh:88a}, \cite{Sh:420}.  So we know that in
 many cases there is set $S^{\bd}_{\lambda,\theta} \subseteq
S^\lambda_\theta := \{\delta < \lambda:\cf(\delta)=\theta\}$ 
such that $\dual(\check I_\theta[\lambda]) = D_\lambda +
 (S^\lambda_\theta \backslash S^{\bd}_{\lambda,\theta})$ and so
$S^{\bd}_{\lambda,\theta}$ is unique (\ref{m13}(4)) modulo
 the club filter, $\cD_\lambda$; for definitions see \S(0C).

We know that consistently, starting with a supercompact we can force
that; e.g. GCH and
$S^{\bd}_{\aleph_{\omega +1},\aleph_n}$ (\ref{m13}(4)) is
stationary for $n=1$ but we do
not know it for $n>1$.  Still this set reflects in no $\aleph_n$, however
we use G.C.H. or just $\aleph_n > 2^{\aleph_0}$.  More generally, 
if $\mu$ is strong limit of cofinality $\aleph_0$
and $S = S^{\bd}_{\mu^+,\aleph_1}$ we do not know if $S$ can reflect
in stationarily many $\delta$'s of cofinality $\aleph_n > \aleph_1$ when
$\aleph_n \le 2^{\aleph_0}$.  Similarly for $\mu$ strong limit of
cofinality $\kappa < \mu$, (see \ref{y1}, \ref{y12}).

By \cite[\S1]{Sh:420} for regular $\lambda,\kappa$ such that $\lambda
> \kappa^+$ there is $S \in \check I_\kappa[\lambda]$ which is
stationary, in fact reflect in stationarily many $\delta < \lambda$ of
cofinality, e.g. $\kappa^{+n} < \lambda$ for $n \ge 1$ (check).  
Related subsets are the
good/bad/chaotic sets of scales ($\langle f_\alpha:\alpha <
\lambda\rangle,f_\alpha \in {}^\kappa \mu$), see \cite[Ch.II]{Sh:g},
\cite{MgSh:204}, \cite{Sh:898} and \ref{m60} here.

The proof in \cite[Ch.IX,\S2]{Sh:g} 
of $\pp(\aleph_\omega) < \aleph_{\omega_4}$ in particular
continue these ideas.  

Recall that if $\bar f = \langle f_\alpha:\alpha < \lambda\rangle$ is 
$<_J$-increasing, $<_J$-cofinal in $\prod\limits_{i < \kappa}
\lambda_i,\lambda_i = \cf(\lambda_i) > \theta \ge \kappa^+$ then
$S^{\ged}_\theta(\bar f) := \{\delta < \lambda:\cf(\delta) = \theta$ and
$\bar f \rest \delta$ is flat (see \ref{m60}) has complement 
orthogonal to $\check I_\theta[\lambda]$ modulo 
the non-stationary ideal, (i.e. have a
non-stationary intersection with any $A \in \check I_\theta[\lambda]$).

We made this work after learning Kojman-Milovich-Spadaro \cite{MK},
which shows $\bar f$ is stable in $\kappa^{+4}$; we learn later that
$S^{\ged}_\theta(\bar f) = S^\lambda_\theta \mod D_\lambda$ when
$\theta = \kappa^{+n},n \ge 4$ was pointed out, following the above, by
Sharon-Viale \cite[footnote 5]{ShVa10}, using Abraham-Magidor
\cite[2.12,2.19]{AbMg10}. 

We start by continuing \cite[\S1]{Sh:420}, \cite[Ch.IX,\S2]{Sh:g},
 to re-examine some of those problems; see \S(0B).  More
specifically, we shed some light on question (Q1) in
\ref{y1}, \ref{y12} proved in \S(1A).

What about (Q2)?  This was a central issue of \cite{Sh:898} using one
dimensional. The $n$-dimensional are from \cite{Sh:883} and lately
\cite{Sh:F1200}, which relies on the results here; 
see \ref{y40}, \ref{y19}, \ref{y40d} and proved in \S(1B).

Much earlier Solovay proved that above a compact cardinal, the
singular cardinal hypothesis holds; it follows that the so called
strong hypothesis $(\mu > \cf(\mu) 
\Rightarrow \pp(\mu) = \mu^+)$ holds; so $\pcf$ becomes
trivial.  Moreover, by \cite[Ch.II]{Sh:g} if 
$\pp_J(\mu) > \lambda = \cf(\lambda) > \mu >
\cf(\mu) = \kappa$ ( where $J \supseteq [\kappa]^{< \kappa}$ is an ideal on
$\kappa$) \then \, there is a sequence $\langle
f_\alpha:\alpha < \lambda\rangle$ with $f_\alpha \in {}^\kappa \mu$
which is $<_J$-increasing and is $\mu^+$-free even as a sequence, 
so $\bar f \rest \delta$ is flat when $\kappa <
\cf(\delta) < \mu$, (i.e. the good set of $\bar f,\ged(\bar f$) is
large.

But if $\kappa = \cf(\mu) < \mu$, the consistency result on
$\check I_{\kappa^+}[\mu^+]$ from \cite{Sh:108} can be strengthened; we know
consistently there are strong reflection properties say if $\GCH$,
consistently the case of Chang conjection holds from $(\aleph_{\omega
  +1},\aleph_\omega) \rightarrow (\aleph_1,\aleph_1)$, by
Levinski-Magidor-Shelah \cite{LMSh:198} and $(\aleph_{\omega + \omega
  +1},\aleph_{\omega + \omega}) \rightarrow (\aleph_{\omega
+1},\aleph_\omega)$.  We can manipulate $2^\kappa$ for $\kappa$ regular.

\subsection {Results} \label{Results} \
\bigskip

What do we accomplish?  E.g. assume $\lambda > \kappa > \aleph_0$ and
for transparency assume $S^{\bd}_{\lambda,\kappa}$ is well defined.  
How much can it reflect?  Assume $\lambda = \mu^+,\cf(\mu) =
\kappa,\mu$ strong limit.  We knew that (\cite{Sh:108}) if, e.g. $\theta =
(2^\kappa)^{+n+1}$ then $S^{\bd}_{\lambda,\kappa}$ 
does not reflect in $S^\lambda_\theta$.
Here \ref{y12} gives more: assuming $(\forall n)(2^{\kappa^{+n}} 
< \lambda)$ we have, e.g. for $n \ge 2,m \ge n+2$: 
if $S^{\bd}_{\lambda,\kappa}$ reflect in 
$S^\lambda_{\kappa^{+n}}$ this
reflection does not reflect in $S^\lambda_{\kappa^{+m}}$; moreover
does not reflect in any $S^\lambda_{\theta^+},\theta \in \Reg \cap
\lambda \backslash \kappa^{+n+2}$.  See more in \ref{y12}.

Returning to e.g. ``if $\bar f$ is $<_J$-increasing cofinal in
$\prod\limits_{i < \kappa} \lambda_i/J$ and $i < \kappa \Rightarrow 
\lambda_i = \cf(\lambda_i) > \kappa$; how large is $S^{\ged}_\theta[\bar f]"$?  
We knew $S^{\ged}_\theta[\bar f]$ is large; here we prove in 
\ref{y1}(1) that:  if
$\theta \in [\kappa^{+4},\kappa^{+\comp(J)}),(\forall i)(\theta <
\lambda_i)$ and $\theta$ is regular $< \lambda$
\then \, $S^{\ged}_\theta[\bar f]$ contains 
$S^\lambda_\theta$ (modulo the club filter of course).   
Hence, e.g. $\bar f$ is $(\theta^{+\comp(J)},\theta^{+4},J)$-free when
$\kappa \le \theta,\theta^{+\comp(J)} < \min\{\lambda_i:i < \kappa\}$,
so if $\lambda_\ell = \pp(\mu_\ell) >
\mu^+_\ell,\mu_\ell > \aleph_0 = \cf(\mu_\ell)$ for $\ell=1,2$ and
$\mu^{+4}_1 \le \lambda_1 < \lambda_2$ \then \, $(\lambda_2,\mu_2) \nrightarrow
(\lambda_1,\mu_1)$.  

\noindent
But this is not enough to prove what we need for Q2, i.e. \ref{y19} which is
$(\theta_2,\theta_1)$-freeness; (the problem being for $\langle
\delta_i:i < \theta\rangle$ increasing continuous, for $i$ of
cofinality $\le \kappa$) but \ref{b40} tells us more, in particular,
enough for Theorem \ref{y19}.

More specifically, we shall show (the proofs are 
given later, the definitions appear in \S(0C) below):
\begin{theorem}
\label{y1}
Assume $\lambda > \sigma > \partial > \theta^+ > \theta > \aleph_0$ 
are regular.

\noindent
1) Some $S \in \check I_\theta[\lambda]$ reflect in every $\delta \in
   S^\lambda_\sigma$, see Definition \ref{m15}(1).

\noindent
2) Moreover, if $\delta \in S^\lambda_\sigma$ \then \, $\{\delta_1 <
\delta:\cf(\delta_1) = \partial$ and $S$ reflects in $\delta_1\}$ 
is a stationary subset of $\delta$.

\noindent
3) Moreover, for any $(\partial,\theta,< \sigma)$-system
   $\bar{\cP}^*$, see Definition \ref{m3}, 
for any ordinal $\delta \in S^\lambda_\sigma$, for
   any increasing continuous sequence $\langle \delta_i:i <
   \sigma\rangle$ of ordinals with limit $\delta$ (clearly exists) for
some $S_1 \in \check I^{\ac}_\partial \langle \sigma,\sigma\rangle$,
   see Definition \ref{m14}(2) we have:
\mn
\begin{enumerate}
\item[$(*)$]  if $j \in S^\sigma_\partial \backslash S_1$ \then \,
  there is $S_2 \in I^{\cg}_\theta(\bar{\cP}^*)$ such that for some
increasing continuous sequence 
$\langle i_\varepsilon:\varepsilon < \partial\rangle$
with limit $j$ we have $\varepsilon \in S^\partial_\theta
\backslash S \Rightarrow \delta_{i_\varepsilon} \in
  \good''_\theta(\bar{\cP})$.
\end{enumerate}
\end{theorem}

\begin{theorem}
\label{y12}
Assume $\lambda > \theta^{+ \omega}$ and $\lambda,\theta$ are regular
uncountable and $2^{\theta^{+n}} < \lambda$ for every $n$.

\noindent
1) If $S^{\bd}_{\lambda,\theta}$ is (well defined and) stationary
\then \, there are $n$ and stationary 
$S \subseteq S^\lambda_{\theta^{+n}}$ which reflects in no ordinal
$\delta$ of cofinality $\in [\theta,\theta^{+\omega})$.

\noindent
2) There is $S \in \check I_\theta[\lambda]$ 
such that for every $n \ge 2$, either $S_1 = 
S^\lambda_{\theta^{+n}} \cap \refl(\lambda \backslash S)$ 
is not stationary (in $\lambda$) \underline{or} $S_1$ is stationary 
but is the union of $\le 2^{\theta^{+n}}$ sets each of 
which reflect in no $\delta$ of cofinality 
$\in [\theta^{n+2},\theta^{+\omega})$. 

\noindent
3) In part (2) in the second possibility some stationary $S_2 \subseteq
   S_1(\subseteq S^\lambda_{\theta^{+n}})$ either reflect in no
   ordinal of cofinality $< \theta^{+ \omega}$ \underline{or} $S_3 =
   \{\delta \in S^\lambda_{\theta^{+n+1}}:S_2 \cap \delta$ is
   stationary in $\delta\}$ is a stationary subset of
   $S^\lambda_{\theta^{+n+1}}$ which reflect in no $\delta < \lambda$ of
 cofinality $< \theta^{+ \omega}$.  
\end{theorem}

\noindent
In \cite{Sh:898} we consider another version of freeness, note that being
$(\theta,\sigma)$-free follows from $\theta$-free and is stronger than
stable in every $\kappa \in [\sigma,\theta)$.  We do not get it fully
  but enough to get ``quite free $\bold k$-combinatorial parameters"
  which is enough for applications in \cite{Sh:F1200},
\cite{Sh:1006}.  

\begin{remark}
\label{y17}
1) Recall that for regular $\partial > \aleph_0,\mu \in \bold C_\partial$
   means $\mu$ is strong limit singular of cofinality $\partial$.

\noindent
2) For $\partial = \aleph_0$ the above is almost equal to (and is
   contained in) the class $\{\mu:\mu > \aleph_0$ strong limit of cofinality
   $\aleph_0\}$, more specifically, the difference does not reflect in
   any singular cardinal.

\noindent
3) Having two possibilities in \ref{y19}, make us prefer the non-tree
   version of the black box, (see \cite{Sh:F1200}).
\end{remark}

\begin{theorem}
\label{y19}
Assume $\sigma < \kappa$ are regular, $\mu \in \bold C_\kappa$,
i.e. $\mu$ is strong limit singular of cofinality $\kappa$.

At last one of the following holds:
\mn
\begin{enumerate}
\item[$(A)$]  there is a $\mu^+$-free $\cF \subseteq {}^\kappa \mu$
of cardinality $\lambda := 2^\mu$,  this is called ``$\mu$ has a 1-solution"
\sn
\item[$(B)$]  $\lambda = 2^\mu$ is regular and there is a 
$(\lambda,\mu,\sigma \times \kappa)-5$-solution, see Definition
\ref{y40}.
\end{enumerate}
\end{theorem}

\begin{claim}
\label{y22}
If $\mu > \kappa = \cf(\mu) > \sigma = \cf(\sigma)$ and we let
$\lambda = \mu^+$ \then \, there is $\bar\eta$ satisfying clauses
(a)-(f) of Definition \ref{y40}.
\end{claim}

\begin{definition}
\label{y40}
Assume $\mu \in \bold C_\kappa,\lambda = 2^\mu = \cf(\lambda),\sigma
= \cf(\sigma) < \kappa$; we say $\bold x$ is a 
$(\lambda,\mu,\kappa,\sigma)-5$-solution \when \, it consists of:
\mn
\begin{enumerate}
\item[$(a)$]  $\bar\eta = \langle \bar \eta_\delta:\delta \in S\rangle$
\sn
\item[$(b)$]  $S \subseteq S^\lambda_\sigma$ is stationary in $\lambda$
  (and $\in \check I_\sigma[\lambda])$
\sn
\item[$(c)$]   $\eta_\delta := \langle \alpha_{\delta,i,j}:(i,j)
  \in \sigma \times \kappa \rangle$ and $\langle \alpha_{\delta,i,0}:i
  < \sigma\rangle$ is increasing with limit $\delta$ and  
$\alpha_{\delta,i,j} \in [\alpha_{\delta,i,0},\alpha_{\delta,i,0} +
  \mu)$ increasing with $j$ and 
$\alpha_{\delta,i,0} + \mu \le \alpha_{\delta,i+1,0}$;
and let $C_\delta = \{\alpha_{\delta,i,j}:(i,j) \in
  \sigma \times \kappa\}$
\sn
\item[$(d)$]  if $\alpha_{\delta_1,i_1,j_1} =
  \alpha_{\delta_2,i_2,j_2}$ then $(i_1,j_1) = (i_2,j_2)$ and $i < i_1
  \wedge j < j_2 \Rightarrow \alpha_{\delta_1,i,j} = \alpha_{\delta_2,i,j}$
\sn
\item[$(e)$]  [freeness] $\bar\eta$ is
  $(\theta^{+\kappa+1},\theta^{+4},J_*)$-free, see \ref{b24}(4) when
  $\kappa \le \theta < \mu$ and $J_* = J^{\bd}_{\sigma \times \kappa}
  = \{u \subseteq \sigma \times \kappa$: for some $(i_*,j_*) \in
  \sigma \times \kappa$ we have $u \subseteq \{(i,j) \in \sigma \times
  \kappa:i < i_*$ and $j<j_*\}$
\sn
\item[$(f)$]  [freeness]  $\bar\eta$ is $(\kappa^+,J_*)$-free
\sn
\item[$(g)$]  [black box]  for every $\chi < \mu$ and $\bar F = \langle
 F_\delta:\delta \in S\rangle$ such that
 $F_\delta:{}^{(C_\delta)}\delta \rightarrow \chi$ there is
$\bar \alpha = \langle \alpha_\delta:\delta \in S \rangle \in {}^S
\chi$ such that $(\forall\eta \in {}^\lambda
\lambda)(\exists^{\stat} \delta \in S)(\alpha_\delta = 
F(\eta \rest C_\delta))$, e.g.
\sn
\item[$(g)'$]  for every relational vocabulary $\tau$ of cardinality
  $< \mu$ there is a sequence $\bar M = \langle M_\delta \in S\rangle,
M_\delta$ a $\tau$-model with universe $C_\delta := \Rang(\eta_\delta)
  = \{\alpha_{\delta,i,j}:i < \sigma,j < \partial\}$ such that for every
$\tau$-model $M$ with universe $\lambda$ we have $(\exists^{\stat} \delta \in
  S)(M_\delta = M \rest C_\delta)$.
\end{enumerate}
\end{definition}

\begin{discussion}
\label{y40d}
1) It may be helpful to use this to prove results by cases.  First, find
a proof using a 1-solution, that is with $\mu^+$-freeness using
(A) of \ref{y19} or at least $\theta_*$-free, 
$\cF \subseteq {}^\kappa \mu,|\cF| = 2^\mu,\theta_*$ large enough so in
\cite{Sh:F1200} terms using $\bold x$ with
$\bold k_{\bold x} = 1$.  Second, use $n$ cases of a $5$-solution (see
\ref{y19}(B) and Definition \ref{y40}) so
have $\bold x = \bold x_0 \times \bold x_1 \times \ldots \times \bold
x_n,\bold x_\ell$ is as above so have enough cases of
$(\theta^\kappa,\theta^{+ 4})$-freeness.  This is done in
\cite{Sh:F1200} which uses Theorem \ref{y19}.

\noindent
2) We may use a different division to cases then \ref{y19}, dividing
   case (B) as in \cite{Sh:898}.  Let $\Upsilon =
   \min\{\partial:2^\partial > 2^\mu\}$; and ask whether $\Upsilon =
   \lambda$ or $\Upsilon < \lambda$.

\noindent
2A) If $\Upsilon = \lambda$ then $\lambda = \lambda^{<\lambda}$ hence
we have better statements on $\lambda$, e.g. if $\lambda$ is a
successor cardinal then we have $\diamondsuit_{S^\lambda_{\aleph_0}}$
or $\diamondsuit_{S^\lambda_{\aleph_1}}$ by \cite{Sh:922}.

\noindent
2B) If $\Upsilon < \lambda$, by \cite[\S2]{Sh:898}, we can construct a
black box for $\Upsilon$ by \cite[\S2]{Sh:898}.
\end{discussion}
\bigskip

\subsection {Quoting Definitions} \label{Quoting} \
\bigskip

\noindent
We try to make this work reasonably self-contained. 

\begin{notation}
\label{m0}
1) For regular uncountable cardinal $\lambda$ let $\cD_\lambda$ be the
   filter generated by the clubs of $\lambda$.

\noindent
2) $\cH(\chi)$ is the set of $x$ with transitive closure of
   cardinality $< \chi$.

\noindent
3) Let $<^*_\chi$ will denote a well ordering of $\cH(\chi)$.

\noindent
4) For regular $\kappa$ and cardinal (or ordinal) $\lambda > \kappa$
   let $S^\lambda_\kappa = \{\delta < \lambda:\cf(\delta) = \kappa\}$.

\noindent
5) For an ideal $J$ on $\kappa$ let $\comp(J)$ be $\sup\{\theta:J$ is
   $\theta$-complete$\}$. 
\end{notation}

\begin{definition}
\label{m3}
1) We say $\bar{\cP}$ is a $(\partial,\theta,< \mu)$-system \when \,:
\mn 
\begin{enumerate}
\item[$(a)$]  $\theta \le \partial$ are regular cardinals
\sn
\item[$(b)$]  $\bar{\cP} = \langle \cP_\alpha:\alpha < \partial \rangle$
\sn
\item[$(c)$]  if $a \in \cP_\alpha$ then $a \subseteq \alpha$ and 
$|a| < \theta$ 
\sn
\item[$(d)$]  $\beta \in a \in \cP_\alpha \Rightarrow a_\alpha \cap
 \beta \in \cP_\beta$
\sn
\item[$(e)$]  $\cP_\alpha$ has cardinality $< \mu$.
\end{enumerate}
\mn
2) If $\mu = \partial$ we may write $(\partial,\theta)$-system.
Instead ``$< \mu^+$" we may write $\mu$.  If $\cP_\alpha =
\{a_\alpha\}$ for $\alpha < \partial$ so $\bar{\cP}$ a $(\partial,<
\theta,1)$-system, and we may write $\bar a = \langle a_\alpha:\alpha <
\partial\rangle$ instead of $\bar{\cP}$.  Instead of $\theta$ we may
write $\le \partial$ when $\theta = \partial^+$.
\end{definition}

\begin{remark}
\label{m4}
Concerning Definition \ref{m3}(1) note that we allow $\mu > \partial$;
in fact, this case was used in \cite[Ch.II]{Sh:g}, in proving: if
$\lambda = \tcf(\prod\limits_{i < \kappa} \lambda_i,<_J),\lambda_i =
\cf(\lambda_i) > \kappa$ and $\mu = \lim_J \langle \lambda_i:i <
\kappa\rangle < \lambda_* =- \cf(\lambda_*) < \lambda$ \then \, there
are $\lambda^*_i = \cf(\lambda^*_i) < \lambda_i$ with $\mu = \lim_J
\langle \lambda^*_i:i < \kappa\rangle$ such that $\lambda_* =
\tcf(\prod\limits_{i < \kappa} \lambda^*_i,<_J)$ exemplified by some
$\mu^+$-free $\langle f_\alpha:\alpha < \lambda_*\rangle$.
\end{remark}

\begin{fact}
\label{m6}
For every regular $\theta$ and stationary $S \subseteq \{\delta <
\theta^+:\cf(\delta) < \theta\}$ there is a 
$(\theta^+,\theta,1)$-system, which means that there is $\bar a$
satisfying:
\mn
\begin{enumerate}
\item[$(a)$]  $\bar a = \langle a_\alpha:\alpha < \theta^+\rangle$
\sn
\item[$(b)$]  $a_\alpha \subseteq \alpha$
\sn
\item[$(c)$]  $|a_\alpha| < \theta$
\sn
\item[$(d)$]  $\beta \in a_\alpha \Rightarrow a_\beta = a_\alpha \cap
\beta$
\sn
\item[$(e)$]  if $E$ is a club of $\theta^+$ and $\zeta < \theta$
\then \, there is $\alpha$ such that $a_\alpha \subseteq E \wedge \alpha =
\sup(a_\alpha) \wedge \otp(a_\alpha) = \zeta$
\sn
\item[$(f)$]  if $E$ is a club of $\theta^+$ and $\zeta < \theta$, 
\then \, for some $\delta \in
S \cap E$ we have $a_\delta \subseteq E \wedge \delta =
\sup(a_\delta)$ and $\zeta$ divides $\otp(a_\delta)$.
\end{enumerate}
\end{fact}

\begin{PROOF}{\ref{m6}}
See \cite[Ch.III]{Sh:g} + correction in \cite{Sh:E12}.
As of guessing clubs for clause (f), it is like \cite[\S1]{Sh:420}.  
We just are more explicit in what we get.
\end{PROOF}

\noindent
Recall (\cite{Sh:108} = \cite{Sh:88a},\cite[\S1]{Sh:420}), (there we
vary $\theta$)
\begin{definition}
\label{m13}
1) Let $\lambda > \theta$ be regular.

\noindent
2) For a $(\lambda,\theta,<\mu)$-system 
$\bar{\cP} = \langle \cP_\alpha:\alpha < \lambda\rangle$ let
\mn
\begin{enumerate}
\item[$\bullet$]  good$'_{<\theta}(\bar{\cP}) = \{\delta <
\lambda:\cf(\delta) < \theta$ and there is an unbounded $u \subseteq
  \delta$ of order type $< \delta$ such that $\alpha \in u \Rightarrow
 u \cap \alpha \in \cP_\alpha\}$
\sn
\item[$\bullet$]  good$''_{<\theta}(\bar{\cP})$ is defined similarly but 
$\otp(u) = \cf(\delta)$.
\end{enumerate}
\mn
2A) For a $(\lambda,\theta,<\mu)$-system $\bar{\cP}$, we define
good$'_{\le \theta}(\bar{\cP})$, good$''_{\le \theta}(\bar{\cP})$
naturally; we defined
$\good'_{=\theta}(\bar{\cP}),\good''_{=\theta}(\bar{\cP})$ similarly
but demand $\cf(\delta) = \theta$ and add ``$u \in \cP_\delta$".

\noindent
3) $\check I_\theta[\lambda]$ is the set of $S \subseteq
S^\lambda_\theta := \{\delta < \lambda:\cf(\delta) = \theta\}$ such
that for some $(\lambda,\theta,1)$-system $\bar a$ 
and club $E$ of $\lambda$ we have $S \cap E \subseteq \text{
  good}'_\theta(\bar{\cP})$, equivalently for some $(\lambda,<
\theta,1)$-system $\bar a$ and club
$E$ of $\lambda,S \cap E \subseteq \text{ good}''_\theta(\bar a)$;
equivalently, we may use $\bar{\cP}$ a $(\lambda,\lambda,<
\lambda)$-system or $(\lambda,\theta,< \lambda)$-system;
abusing notation for $S \subseteq \lambda,S \in \check
I_\theta[\lambda]$ means $S \cap S^\lambda_\theta \in \check
I_\theta[\lambda]$.

\noindent
4) If $\check I_\theta[\lambda] =$ (the non-stationary ideal on
$S^\lambda_\theta$) $+ S_*$ \then \, we call $S_*$ the
good set on $\lambda$ for cofinality $\theta$; it will be denoted
$S^{\ged}_{\lambda,\theta}$; its complement $S^{\bd}_{\lambda,\theta}
:= S^\lambda_\theta \backslash S_*$ is called the bad set; of course,
as only $S_*/\cD_\lambda$ is 
unique this notation pedentically is not justified.

\noindent
5) Let $\check I^\perp_\kappa[\lambda] = \{S \subseteq S^\lambda_\kappa$:
   if $S_1 \in \check I_\kappa[\lambda]$ then $S_1 \cap S$ is not
   stationary (in $\lambda$)$\}$.

\noindent
6) Let $\check I[\lambda] = \{S \subseteq \lambda$: if $\theta =
\cf(\theta) < \lambda$ then $S \cap S^\lambda_\theta \in \check
I_\theta[\lambda]\}$. 
\end{definition}

\begin{definition}
\label{m14}
Let $\lambda > \theta^+$ be regular. 

\noindent
1) Let $I^{\cg}_{\theta}[\lambda,\mu]$ be the set of $S \subseteq
   S^\lambda_\theta$ such that ($\cg$ stands for club guessing) there
is no $\bar{\cP}$ witnessing $S \in
   (I^{\cg}_\theta[\lambda,\mu])^+$ which means $S \subseteq
   S^\lambda_\theta \wedge S \notin I^{\cg}_\theta(\bar{\cP})$ that is:
\mn
\begin{enumerate} 
\item[$(*)_1$]  $\bar{\cP} = \langle \cP_\alpha:\alpha <
  \lambda\rangle$ is a $(\lambda,\theta,< \mu)$-system
\sn
\item[$(*)_2$]  for $\bar{\cP},\lambda$ as above let
$I^{\cg}_\theta(\bar{\cP})$ be the set of $S \subseteq \lambda$ such that
\sn
\begin{enumerate}
\item[$\bullet$]  for some club $E$ of $\lambda$ for no
  $\delta \in S$ and $a \in \cP_\delta$ do we have $a \subseteq E \wedge
  \sup(a) = \delta$.
\end{enumerate}
\end{enumerate}
\mn
1A) We define $I^{\dg}_\theta[\lambda,\mu],I^{\dg}_\theta(\bar{\cP})$
similarly except that in $\bullet$ of $(*)_2$ we demand only 
$a \in \cP_{< \lambda}$.

\noindent
2) Assume $\lambda = \cf(\lambda) \ge \theta = \cf(\theta),\lambda \ge
\mu,\mu^+ \ge \theta$.  Let $\check I^{\ac}_\theta\langle \lambda,
\mu\rangle$ be the set of $S \subseteq S^\lambda_\theta$ such that
 there are $\chi > \lambda + \mu$ and $x \in \cH(\chi)$ for which 
there is no sequence $\bar N = 
\langle N_\varepsilon:\varepsilon < \theta\rangle$ satisfying:
\mn
\begin{enumerate}
\item[$(a)$]  $N_\varepsilon \prec (\cH(\chi),\theta,<^*_\chi)$
\sn
\item[$(b)$]  $\langle N_\zeta:\zeta < \theta\rangle$ is increasing
  continuous
\sn
\item[$(c)$]  $\langle N_\zeta:\zeta \le \varepsilon \rangle \in
  N_{\varepsilon +1}$
\sn
\item[$(d)$]  $\|N_\varepsilon\| < \mu$ and $N_\varepsilon \cap \mu$
  is an ordinal
\sn
\item[$(e)$]  $\{x,\lambda,\mu,\theta\} \in N_0$
\sn
\item[$(f)$]  $\cup\{N_\varepsilon \cap \lambda:\varepsilon < \theta\} \in
  S$.
\end{enumerate}
\end{definition}

\begin{definition}
\label{m15}
For $\lambda$ regular uncountable and unbounded $S \subseteq \lambda$
let $\refl(S) = \{\delta < \lambda:\cf(\delta) > \aleph_0$ and $S$
reflects in $\delta\}$ where ``$S$ reflects in $\delta$" means $S \cap
\delta$ is a stationary subset of $\delta\}$.

\noindent
2) We say $S \subseteq \lambda$ reflects in $S^\lambda_\theta$ if
   $\{\delta \in S^\lambda_\theta:S \cap \delta$ is stationary in
   $\delta\}$ is a stationary subset of $\lambda$.  We may replace
   $S^\lambda_\theta$ by any stationary subset of $\lambda$. 
\end{definition}

\begin{definition}
\label{m17}
For a regular cardinal $\partial$, let $\bold C_\partial$ be the class
of strong limit singular cardinals $\mu$ of cofinality $\partial$ such that
$\pp^*(\mu) =^+ 2^\mu$.
\end{definition}

\begin{fact}
\label{m19}
1) Assume $\lambda$ is regular and $\lambda = \cf(\lambda) > \mu$ and
   if $\lambda = \mu^+,\mu$ regular, \underline{or} 
$\alpha < \lambda \Rightarrow \cf([\alpha]^{< \mu},\subseteq) <
   \lambda$ and $\mu^+ < \lambda$, \then \, $\theta = \cf(\theta) <
   \mu \Rightarrow S^\lambda_\theta \in \check I_\theta[\lambda]$,
   moreover, there is a closed $(\lambda,\mu,< \lambda)$-system 
$\bar{\cP}$ such that:  $\delta < 
\lambda \wedge \cf(\delta) < \mu \Rightarrow (\exists a \in \cP_\delta)
(\sup(a) = \delta \wedge \otp(a) = \cf(\delta))$.

\noindent
2) $\check I^{\ac}_\theta\langle \lambda,\mu \rangle \subseteq \check
   I_\theta[\lambda]$ when defined.

\noindent
3) If $\lambda > \theta^+$ and $\lambda,\theta$ are regular \then \,
   there is a $(\lambda,\le \theta,< \lambda)$-system $\bar{\cP}$ such
   that $S^\lambda_\theta \notin I^{\cg}_\theta(\bar{\cP})$ and
$\otp(a) = \theta$.
\end{fact}

\begin{discussion}
\label{m52}
1) For the equivalence of the two versions in Definition \ref{m13}(3),
   see \cite[\S1]{Sh:420}.

\noindent
2) When does $S^{\ged}_{\lambda,\theta}$ exist?

See \cite{Sh:108} = \cite{Sh:88a}, $S^{\ged}_{\lambda,\theta}$ exists
under quite weak cardinal arithmetic assumptions (much weaker than GCH).

\noindent
3) Of course, if $\alpha < \lambda \Rightarrow |\alpha|^{< \theta} <
 \lambda$ then $S^{\bd}_{\lambda,\theta} = \emptyset$.

\noindent
4) It is proved there for $\lambda$, e.g. successor of strong limit
 singular $\mu$ and $\theta \in (\cf(\mu),\mu)$ that
 $S^{\bd}_{\lambda,\theta}$ exists and
does not reflect in cofinality $(2^\theta)^+$ and in cofinality
$\partial$ when $(\forall \alpha < \partial)[|\alpha|^\theta < \partial]$.

\noindent
5) Also it is proved (\cite[Ch.II]{Sh:g}) that if $\lambda$ is a 
successor of regular $\aleph_0 < \theta = \cf(\theta)$ and 
$\theta^+ < \lambda$ then  $S^{\bd}_{\lambda,\theta}$ is $\emptyset$; 
(i.e. not stationary), see \ref{m6} or \ref{m19}(1).
\end{discussion}

\noindent
In \S(1B) we shall use \cite[Ch.II]{Sh:g}.
\begin{definition}
\label{m60}
Let $\bar f$ be $<_J$-increasing in ${}^\kappa\Ord,J$ an ideal on $I$.

\noindent
1) We say $\bar f$ is flat in $\delta$ or $\delta \in S_{\ged}[\bar
  f,J] = S^{\ged}_J[\bar f]$ 
\when \, $\delta \le \ell g(\bar f),\cf(\delta) > \kappa$ and 
there is a $<_J$-eub $g$ to $\bar f \rest \delta$ such that 
$(\forall i < \kappa)(\cf(g(i)) = \cf(\delta))$, 
equivalently there are increasing sequences $\langle
\alpha_{i,\varepsilon}:\varepsilon < \cf(\delta)\rangle$ for $i <
   \kappa$ such that $(\forall \alpha < \delta)(\exists \varepsilon <
\cf(\delta))(f_\alpha <_J \langle \alpha_{i,\varepsilon}:i <
   \kappa\rangle)$ and $(\forall \varepsilon < \cf(\delta))(\exists
   \alpha < \delta)(\langle \alpha_{i,\varepsilon}:i < \kappa\rangle
   <_J f_\alpha)$.

\noindent
2) We say $\delta$ is strongly chaotic for $\bar f$ or $\delta \in
S_{\sch}[\bar f,J] = S^{\sch}_J[\bar f]$ \when \,
there is a sequence $\langle u_i:i < \kappa\rangle,
u_i \subseteq \Ord,|u_i| \le \kappa$
 and $(\forall \alpha < \delta)(\exists g \in \prod\limits_{i}
  u_i)(\exists \beta < \delta)(f_\alpha <_J g <_J f_\beta)$.

\noindent
2A) We say $\delta$ is chaotic for $\bar f$ or $\delta \in
S^{\ch}_J[f] = S_{\ch}[\bar f,J]$ \when \, there is $\bar u$ as above
such that for every $\alpha < \delta$ for some $\beta \in
(\alpha,\delta)$ the set $A_{\alpha,\beta} = A_{\alpha,\beta}[\bar
  u,\bar f]$ belongs to $J^+$ where $A_{\alpha,\beta} =
\{i < \kappa:\min(u_i \cup \{\infty\}
\backslash f_\alpha(i)\} < \min(u_i \cup \{\infty\} \backslash f_\beta(i))\}$.

\noindent
2B) We define $S^{\sch}_\theta[\bar f,J] = S^{\sch}_{J,\theta}[\bar
  f],S^{\ch}_\theta[\bar f,J] = S^{\ch}_{J,\theta}[\bar f]$ similarly
  but restricting ourselves to $\delta$ of cofinality $\theta$.

\noindent
3) We say $\delta$ is bad for $\bar f$ or $\delta \in S_{\bd}
[\bar f,J] = S^{\bd}_J[\bar f]$ \when \, $\delta \le \ell g(\bar
f),\cf(\delta) > \kappa$ and $\bar f \rest \delta$ has 
$<_J$-eub $g$ but is not flat.
\end{definition}

\begin{claim}
\label{m62}
Let $J,\bar f$ be as in \ref{m60}.

\noindent
1) If $\delta \le \ell g(\bar f),\cf(\delta) > \kappa^+$ \then \,
   $\delta$ satisfies exactly one of good, bad or chaotic.

\noindent
2) In other words $\{\delta:\delta \le \ell g(\bar f)$ and $\cf(\delta) >
\kappa^+\}$ is included in the disjoint union of $S_{\ged}[\bar f],
S_{\bd}[\bar f],S_{\sch}[\bar f]$.
\end{claim}

\begin{PROOF}{\ref{g5}}
By \cite[Ch.II,\S2]{Sh:g}.
\end{PROOF}

\begin{claim}
\label{m65}
Let $\bar f,J,\kappa$ be as in \ref{m60} and $\lambda = \ell g(\bar f)$.

\noindent
1) If $\delta \in S^{\ch}_J[\bar f]$ \then \, for some club $e$ of $\delta$,
we have $\alpha \in e \wedge \cf(\alpha) > \kappa \Rightarrow \alpha 
\in S^{\ch}_J[\bar f]$.

\noindent
1A)  Similarly for $S_{\sch}[\bar f]$.

\noindent
2) If $\delta \in S^{\ged}_J[\bar f]$ \then \, for some club $e$ of $\delta$
we have $\alpha \in e \wedge \cf(\alpha) > \kappa \Rightarrow
\alpha \in S^{\ged}_J[\bar f]$.

\noindent
3) If $\delta \le \lambda,\cf(\delta) \in S_{\bd}[\bar f]$ then
   $\cf(\delta) \ge \kappa^{+ \comp(J) +1}$.
\end{claim}

\begin{claim}
\label{m68}
Assume $(\lambda,\bar\lambda,J,\kappa)$ is a $\pcf$ case, $\bar f$ a
witness for it, see Definition \ref{b18}.  
If $\kappa < \sigma < \min\{\lambda_i:i < \kappa\}$ or just
$\kappa < \sigma < \lim-\inf_J(\bar\lambda)$ and $S \in \check
I_\sigma[\lambda]$ \then \, $E \cap S \subseteq S_{\ged}[\bar f]$ for some
club $E$ of $\lambda$.
\end{claim}
\newpage

\section {On systems} \label{On}
\bigskip

\subsection {Existence of large members of $\check I_\theta[\lambda]$}
\label{Existence} \
\bigskip

\begin{claim}
\label{b4}
Assume $\lambda > \aleph_1$ is regular and $M_* \prec (\cH(\lambda),\in)$ has
cardinality $< \lambda$ and $\{\lambda,\theta\} \subseteq M_*$ and $M_*
\cap \lambda \in \lambda$.
\Then \, we can find a pair $(E,\bar{\cP})$ which is
$(\lambda,M_*)$-suitable, which means:
\mn
\begin{enumerate}
\item[$\boxplus$]  $(a) \quad E$ is a club of $\lambda$; we may add
  $\alpha \in E \wedge \alpha > \sup(\alpha \cap E) \Rightarrow
  \cf(\alpha) = \aleph_0$
\sn
\item[${{}}$]  $(b) \quad \bar{\cP} = \langle \cP_\alpha:\alpha <
  \lambda\rangle$ is a $(\lambda,\lambda,< \lambda)$-system and $\theta =
\cf(\theta) < \lambda \cap M_* \Rightarrow$

\hskip25pt $\good''_\theta(\bar{\cP}) \supseteq
  S^\lambda_\theta \backslash E$ 
\sn
\item[${{}}$]  $(c) \quad$ if $\sigma > \partial$ are
regular $\in \lambda \cap M_*$ and

\hskip25pt $\bar{\cP}^* = 
\langle \cP^*_\alpha:\alpha < \partial\rangle \in M_*$ is a 
$(\partial,\partial,< \sigma)$-system and

\hskip25pt $\langle \delta_i:i < \sigma \rangle$
is an increasing continuous sequence of members 

\hskip25pt  from $E$, \then \, there are $f,e$ such that:
\sn
\begin{enumerate}
\item[${{}}$]  $(\alpha) \quad e$ is a club of $\partial$
\sn
\item[${{}}$]  $(\beta) \quad f$ is an increasing continuous function
from $\partial$ into $\{\delta_i:i < \sigma\}$
\sn
\item[${{}}$]  $(\gamma) \quad$ if $\varepsilon < \partial,a \in
\cP^*_\varepsilon$ and $a \subseteq e$ \then \, 
$\{f(\xi):\xi \in a$ and $\otp(a \cap \xi)$

\hskip25pt is a successor ordinal$\} \in \cP_{f(\varepsilon +1)}$
\end{enumerate}
\sn
\item[${{}}$]  $(c)^+ \quad$ like (c) but we replace $(\gamma)$ by
\sn
\begin{enumerate}
\item[${{}}$]  $(\gamma)^+ \quad$ if $\varepsilon < \partial,
a \in \cP^*_\varepsilon$ and
$a \subseteq e$ and $\langle \gamma_\iota:\iota < \otp(a)\rangle$ list 
$a$

\hskip30pt   in increasing order \then \, in addition 
to the conclusion of $(\gamma)$
\sn
\item[${{}}$]  \hskip10pt $\bullet \quad$ we can choose $\beta_\iota \in 
[\gamma_\iota,\gamma_{\iota +1})$ for $\iota < \otp(a)$ such that

\hskip25pt  $\{\beta_j:j \le \iota\} \in \cP_{\beta_{\iota +1}}$ for 
every $\iota < \otp(a)$
\sn
\item[${{}}$]  \hskip10pt $\bullet \quad$ if $a$ has no last member
  then $\sup(a) \in \good''_\theta(\bar{\cP})$
\end{enumerate}
\sn
\item[${{}}$]  $(d) \quad$ if $\langle \delta_i:i < \sigma\rangle$ is
  an increasing continuous sequence of members of $E$

\hskip25pt  and $\sigma > \partial > \theta$ are regular $\in \lambda
\cap M_*$ and $\bar{\cP}^* = \langle \cP^*_\varepsilon:
\varepsilon < \partial\rangle \in M_*$ 

\hskip25pt  is a $(\partial,\le \theta,< \sigma)$-system
\then \, for some $e,f$ satisfying 

\hskip25pt clauses $(\alpha),(\beta),(\gamma),(\gamma)^+$ we have
\sn
\begin{enumerate}
\item[${{}}$]  $(\delta) \quad$ the following set belongs to
  $I^{\dg}_\theta(\bar{\cP}^*)$, recalling \ref{m14}(1A)

\hskip25pt  $\{\zeta \in S^\partial_\theta$: there
  is no $a \subseteq e,a \in \cP^*_{< \partial}$ such that $a
  \subseteq \zeta = \sup(a)$

\hskip25pt and $\otp(a) = \theta\}$
\sn
\item[${{}}$]  $(\varepsilon) \quad$ the following set belongs to
  $\check I^{\ac}_\partial \langle \sigma,\sigma\rangle$, see
  Definition \ref{m14}(2)

\hskip25pt  $\{i \in S^\sigma_\partial$: there are no $e,f$ satisfying
  $\sup(e)=i$ and

\hskip25pt clauses $(\alpha),(\beta),(\gamma),(\gamma)^+,(\delta)$ above$\}$. 
\end{enumerate}
\end{enumerate}
\end{claim}

\begin{remark}
\label{b5}
1) Note that for $\good''_\theta(\bar{\cP})$, only $\langle \cP_\alpha
   \cap [\alpha]^{< \theta}:\alpha < \lambda\rangle$ matters.

\noindent
2) For $\bar M$ as in $\odot_1$ in the proof and $\alpha < \lambda$
essentially $\bar{\cP}$ satisfies the conclusion with $M_*$ replaced
by $M_\alpha$, the essentially because we should ignore the ordinals $\le
\alpha$, i.e. in clauses $(c),(c)^+,(d)$ demand $\delta_0 > \alpha$.
\end{remark}

\begin{PROOF}{\ref{b4}}
Let $\chi > \lambda$ and let $\bar M$ be such that:
\mn
\begin{enumerate}
\item[$\odot_1$]  $(a) \quad \bar M = 
\langle M_\alpha:\alpha < \lambda\rangle$ be a 
$\prec$-increasing continuous sequence
\sn
\item[${{}}$]  $(b) \quad M_\alpha \prec (\cH(\chi),\in,<^*_\chi)$
\sn
\item[${{}}$]  $(c) \quad \|M_\alpha\| < \lambda$
\sn
\item[${{}}$]  $(d) \quad \bar M \rest (\alpha +1) \in 
M_\alpha$
\sn
\item[${{}}$]  $(e) \quad M_\alpha \cap \lambda \in
\lambda$ for every $\alpha < \lambda$
\sn
\item[${{}}$]  $(f) \quad \bar{\cP}^* \in M_0$
\sn
\item[${{}}$]  $(g) \quad$ moreover $M_* \in M_0$ hence $M_* \subseteq M_0$.
\end{enumerate}
\mn
Let $E = \{\alpha:M_\alpha \cap \lambda = \alpha\}$.  Clearly $E$ is a club
of $\lambda$, hence clause (a) of $\boxplus$ holds, even the ``we may add".

Let $\bar{\cP} = \langle \cP_\alpha:\alpha < \lambda\rangle$ be
defined by:
\mn
\begin{enumerate}
\item[$\odot_2$]   $\cP_\alpha = \{a \in M_{\alpha +1}:a \subseteq
  \alpha$ so $|a| < \lambda$ and $\beta \in a \Rightarrow a \cap \beta \in
M_{\beta +1}\}$ so $\bar{\cP} =
\langle \cP_\alpha:\alpha < \lambda\rangle$ is a
$(\lambda,\lambda,<\lambda)$-system, moreover, $\boxplus(b)$ holds.
\end{enumerate}
\mn
[Why does $\boxplus(b)$ hold?  Let $\delta \in 
S^\lambda_\theta \backslash E$ be a limit
ordinal, so for some $\alpha < \delta$ we have $\delta \in M_\alpha$
hence there is an unbounded (and even closed) subset $a$ of $\delta$
  in $M_\alpha$ of order type $\cf(\delta)$ 
so $\beta \in (a \backslash \alpha) \Rightarrow (a 
\backslash \alpha) \cap \beta \in M_\alpha \subseteq M_\beta
\Rightarrow (a \backslash \alpha) \cap \beta \in M_\beta$.  So
  indeed $\good''_\theta(\bar{\cP}) \supseteq S^\lambda_\theta \backslash E$.]

So we arrive to the main point, that is to prove clauses $(c),(c)^+$
and later comment on its relative $(d)$.
So let $\partial < \sigma \in M_* \cap \lambda$ be regular and
$\bar{\cP}^* \in M_*$ be a $(\partial,\partial,< \sigma)$-system and
let $\bar\delta = \langle \delta_i:i < \sigma \rangle$ be an
increasing continuous sequence of ordinals from $E$ and let
$\delta_\sigma := \cup\{\delta_i:i < \sigma\}$ so also $\langle
\delta_i:i \le \sigma\rangle$ is an increasing continuous sequence of
ordinals from $E$.

We choose $N_\varepsilon$ by induction on $\varepsilon \le \partial$
such that:
\mn
\begin{enumerate}
\item[$\odot_3$]  $(a) \quad N_\varepsilon \prec
(\cH(\chi),\in,<^*_\chi)$
\sn
\item[${{}}$]  $(b) \quad \|N_\varepsilon\| < \sigma$
\sn
\item[${{}}$]  $(c) \quad \langle N_\xi:\xi \le \zeta\rangle \in
N_\varepsilon$ when $\zeta < \varepsilon$
\sn
\item[${{}}$]  $(d) \quad \langle N_\zeta:\zeta \le
\varepsilon\rangle$ is $\prec$-increasing continuous
\sn
\item[${{}}$]  $(e) \quad \lambda,\sigma,\partial,\theta,
E,\bar M,\bar\delta$ and $\cP^*$ belongs to
$N_\varepsilon$
\sn
\item[${{}}$]  $(f) \quad \partial +1 \subseteq N_\varepsilon$
  moreover (follows if $\sigma = \partial^+$)
$N_\varepsilon \cap \sigma \in (\partial,\sigma)$.
\end{enumerate}
\mn
This is easy.  Let $i(\varepsilon) := N_\varepsilon \cap \sigma$
for $\varepsilon \le \partial$, hence $i(\varepsilon) < \sigma$ is
increasing continuous with $\varepsilon$.  So
$\delta_{i(\varepsilon)}$ is an ordinal $\in E \subseteq \lambda$ hence
$M_{\delta_{i(\varepsilon)}}$ is well defined and
$\delta_{i(\varepsilon)} \in M_{\delta_{i(\varepsilon)}+1}$, also
$\langle \delta_{i(\varepsilon)}:\varepsilon < \partial\rangle$ is
increasing continuous with limit $\delta_{i(\partial)}$.  For
$\varepsilon = \partial$ clearly $\cf(\delta_{i(\varepsilon)}) =
\cf(\delta_{i(\partial)}) = \cf(\partial) = \partial$ hence
\mn
\begin{enumerate}
\item[$\oplus_1$]  $(a) \quad$ there 
is a club $C$ of $\delta_{i(\partial)}$ of order type 
$\cf(\delta_{i(\partial)}) = \partial$
\sn
\item[${{}}$]  $(b) \quad$ necessarily $C \in
\cH(\chi)$ and \wilog \, $C \in M_{\delta_{i(\partial)}+1}$
\sn
\item[${{}}$]  $(c) \quad$ let $g$ be the unique 
increasing continuous function from 
$\partial$ onto $C$, so 

\hskip25pt  necessarily $g \in M_{\delta_{i(\partial)}+1}$
\sn
\item[${{}}$]  $(d) \quad$ let $e = 
\{\varepsilon < \partial:\delta_{i(\varepsilon)} \in C$,
moreover $\varepsilon = \otp(C \cap \delta_{i(\varepsilon)})$ and,

\hskip25pt  actually follows, $\delta_{i(\varepsilon)} = g(\varepsilon)\}$
\sn
\item[${{}}$]  $(e) \quad$ let $f:\partial \rightarrow \sigma$ be defined by
$f(\varepsilon) = \delta_{i(\varepsilon)}$.
\end{enumerate}
\mn
Now $C$ is a club of $\partial$ and both 
$\langle g(\varepsilon):\varepsilon < \partial \rangle$ and
$\langle \delta_{i(\varepsilon)}:\varepsilon < \partial\rangle$ are
increasing continuous sequences of ordinals with limit
$\delta_{i(\partial)}$, so clearly
\mn
\begin{enumerate}
\item[$\oplus_2$]   $e$ is a club of $\partial$.  
\end{enumerate}
\mn
So concerning clause (c) (of $\boxplus$) it suffices to 
prove that the pair $(f,e)$ we
have just chosen is as required there.  Now obviously
$e,f$ satisfy sub-clauses $(\alpha),(\beta)$ of (c).  What
about sub-clause $(\gamma)$ of clause (c) and subclause $(\gamma)^+$ 
of clause (c)$^+$?  

Clearly
\mn
\begin{enumerate}
\item[$\oplus_3$]   $f \rest e = g \rest e$, see the definition of $e$.
\end{enumerate}
\mn
Now we shall prove
\mn
\begin{enumerate}
\item[$\oplus_4$]   if $\varepsilon < \partial$ and $a \in \cP^*_\varepsilon$ 
satisfies $a \subseteq e$, \then \, $\{g(\zeta):\zeta \in a\} \in
M_{f(\varepsilon +1)}$.
\end{enumerate}
\mn
The proof of $\oplus_4$ is done in $(*)_{4.1} - (*)_{4.7}$.

Note
\mn
\begin{enumerate}
\item[$(*)_{4.1}$]   $\cP^*_\varepsilon \subseteq N_0 \cap M_0
 \subseteq N_{\varepsilon +1} \cap M_{\delta(\partial)+1} \subseteq
 N_{\varepsilon +1} \cap M_{\delta_\sigma}$.
\end{enumerate}
\mn
[Why?  Obvious as $\bar{\cP}^* \in M_*,\partial = \ell g(\bar{\cP}^*) 
\in M_* \cap \lambda$ but $M_* \cap \lambda \subseteq 
M_0 \cap \lambda \in \lambda$ hence $\partial \subseteq M_0$ 
so together $\cP^*_\varepsilon \in M_0$.  
Now $|\cP^*_\varepsilon| < \sigma < \lambda$ and $\sigma \in M_* \cap
\lambda \subseteq M_0 \cap \lambda \in \lambda$ 
so $\cP^*_\varepsilon \subseteq M_0 \subseteq
M_{\delta_{i(\varepsilon)}} \subseteq M_{i(\partial)} \subseteq
M_{\delta_\sigma}$.  Also $\bar{\cP}^* \in N_0$ and $\varepsilon,\partial
\in N_{\varepsilon +1}$ and $|\cP^*_\varepsilon| + \partial < \sigma$
and by $\odot_3(f)$ we have
$N_{\varepsilon +1} \cap \sigma \in \sigma$
hence $\cP^*_\varepsilon \subseteq N_{\varepsilon +1}$,
so together we are done.]   

Also
\mn
\begin{enumerate}
\item[$(*)_{4.2}$]  $\{g(\zeta):\zeta \in a\} \in
M_{\delta_{i(\partial)+1}} \prec M_{\delta_\sigma}$.
\end{enumerate}
\mn
[Why?  As $a$ and $g$ belong to this model; why?  For $a$ because $a
  \in \cP^*_\varepsilon$, see the assumption of $\oplus_4$ and 
$\cP^*_\varepsilon \subseteq M_0 \subseteq M_{\delta_{i(\partial)}} \subseteq
  M_{\delta_{i(\partial)+1}}$ by $(*)_{4.1}$.  For $g$, by the choice
  of $C$ and $g$, see $\oplus_1(a),(b),(c)$.]
\mn
\begin{enumerate}
\item[$(*)_{4.3}$]   $\{g(\zeta):\zeta \in a\} = 
\{(f \rest \varepsilon)(\zeta):\zeta \in a\} \in N_{\varepsilon +1}$.
\end{enumerate}
\mn
[Why?  The equality holds by $\oplus_3$ as $a \subseteq e \wedge a
  \subseteq \varepsilon$ by the
  assumption of $\oplus_4$.  The membership ``$\in N_{\varepsilon
    +1}$" holds as on the one hand $a \subseteq \varepsilon,a \in
  \cP^*_\varepsilon$ hence by $(*)_{4.1}$ also 
$a \in N_{\varepsilon +1}$ and on the other hand
$f \rest \varepsilon \in N_{\varepsilon +1} \prec N_\partial$ as 
$\langle N_\zeta:\zeta \le \varepsilon \rangle \in
N_{\varepsilon +1}$ by $\odot_3(c)$ hence $\langle i(\zeta):\zeta \le
\varepsilon\rangle \in N_{\varepsilon +1}$ by the choice $i(\zeta) =
\sup(N_\zeta \cap \sigma)$ after $\odot_3$ and $\bar\delta \in N_0$ by
$\odot_3(e)$ hence $\langle \delta_{i(\zeta)}:\zeta \le
\varepsilon\rangle \in N_{\varepsilon +1}$ so $f \rest (\varepsilon
+1) \in N_{\varepsilon +1}$ by $\oplus_2(e)$.]

As $\bar\delta \in N_0 \prec N_{i(\partial)}$ by $\odot_3(e)$ we have
$\bar\delta =  \langle \delta_i:i \le \sigma\rangle \in N_0
\prec N_{\varepsilon +1}$ so necessarily $\delta_\sigma
 \in N_0 \prec N_{\varepsilon +1}$ and recalling $\bar M \in N_0$ by
 $\odot_3(e)$ it follows that $M_{\delta_\sigma} = \cup\{M_\alpha:\alpha <
 \delta_\sigma\} \in N_{\varepsilon +1}$ and 
$\bar M \rest \delta_\sigma \in N_{\varepsilon +1} \prec
(\cH(\chi),\in,<^*_\chi)$ hence
\mn
\begin{enumerate}
\item[$(*)_{4.4}$]   $M_{\delta_\sigma} \cap N_{\varepsilon +1}
\subseteq M_{\sup(N_{\varepsilon +1} \cap \delta_{i(\partial)})}$
\end{enumerate}
\mn
but (by $(*)_{4.2} + (*)_{4.3}$) 
\mn
\begin{enumerate}
\item[$(*)_{4.5}$]   $\{g(\zeta):\zeta \in a\} \in M_{\delta_\sigma}
\cap N_{\varepsilon +1}$.
\end{enumerate}
\mn
Now as $\bar M,\bar\delta \in N_0$ and $\sigma \in N_0$ by $\odot_3(e)$, 
clearly $M_{\delta_\sigma} \in N_0$ and
as $N_{\varepsilon +1} \cap \sigma = i(\varepsilon +1)$ by the choice
of $i(\varepsilon +1)$ after $\odot_3$ and
$\|N_{\varepsilon +1}\| < \sigma$ by $\odot_3(b)$ clearly
\mn
\begin{enumerate}
\item[$(*)_{4.6}$]  $N_{\varepsilon +1} \cap M_{\delta_\sigma} \subseteq
M_{\delta_{i(\varepsilon +1)}}$.
\end{enumerate}
\mn
But $f(\varepsilon +1) = \delta_{i(\varepsilon +1)}$ by $\oplus_1(e)$
 hence by $(*)_{4.5} + (*)_{4.6}$ we have
\mn
\begin{enumerate}
\item[$(*)_{4.7}$]  $\{g(\zeta):\zeta \in a\} \in M_{f(\varepsilon +1)}$.
\end{enumerate}
\mn
So we have proved $\oplus_4$.
\mn
\begin{enumerate}
\item[$\oplus_5$]  if $\varepsilon < \partial,a \in
  \cP^*_\varepsilon,a \subseteq e$ and $\xi \in a \wedge (a \cap
\xi$ has a last member) then $\{g(\zeta):\zeta \in a \cap \xi\} \in
  M_{f(\xi)}$.
\end{enumerate}
\mn
[Why?  Let $\zeta(*) = \max(a \cap \xi)$, it is well defined by the
assumption on $\xi$.  But $\bar{\cP}^*$ is a $(\partial,\partial,<
\sigma)$-system by the assumption of clause (c) (so of clause $(c)^+$)
of $\boxplus$, hence by clause (d) of Definition \ref{m3}(1) we have
$a \cap \zeta(*) \in \cP^*_{\zeta(*)}$ and,
  of course, $a \cap \zeta(*) \subseteq e$ hence we can apply
  $\oplus_4$ with $(\zeta(*),a \cap \zeta(*))$ here standing for
  $(\varepsilon,a)$ there, so we can deduce $\{g(\zeta):\zeta \in a
  \cap \zeta(*)\} \in M_{f(\zeta(*)+1)}$.  But $\zeta(*) +1 \le \xi$
hence $f(\zeta(*)+1) \le f(\xi)$ hence $M_{f(\zeta(*)+1)} \subseteq
  M_{f(\xi)}$.  So $\{g(\zeta):\zeta \in a \cap \zeta(*)\} \in
  M_{f(\xi)}$, hence by the obvious closure properties of $M_{f(\xi)}
  \cap [f(\xi)]^{\le \theta}$ also $\{g(\zeta):\zeta \in a \cap \xi\}
  \in M_{f(\xi)}$.]
\mn
\begin{enumerate}
\item[$\oplus_6$]  if $\varepsilon < \partial,a \in \cP^*_\varepsilon$
  and $a \subseteq e$ \then \, the set $b = \{f(\zeta):\zeta \in a$
  and $\otp(a \cap \zeta)$ is a successor ordinal$\}$ belongs to
  $\cP_{f(\varepsilon +1)}$.
\end{enumerate}
\mn
[Why?  By $\oplus_4 + \oplus_5$, the definition of $\cP_{f(\varepsilon
 +1)}$ in $\odot_2$ and the obvious closure properties of each $M_\alpha$.] 

So we are done proving clause $(c)(\gamma)$ of $\boxplus$ hence clause
(c).  Clause $(c)^+ (\gamma)^+$ is proved similarly.

We are left with proving clause (d) of $\boxplus$, let $x =
\{\lambda,\sigma,\partial,\theta,\bar{\cP}^*,E,\bar M\}$ and let $S_1 = \{j
\in S^\sigma_\partial$: there is $\bar N$ as in $\odot_3$
such that $j = \sup(\cup\{N_\varepsilon:\varepsilon < \partial\} \cap
\sigma)\}$.  Now by the definition \ref{m14}(2) of $\check
I^{\ac}_\partial \langle \sigma,\sigma\rangle$ we know that
$S^\sigma_\theta \backslash S_1 \in \check I^{\ac}_\partial \langle
\sigma,\sigma \rangle$.

Next, for each $j \in S_1$ let $\langle N_\varepsilon:\varepsilon
< \partial\rangle$ witness that $j \in S_1$.  Now choose $C,g,e,f$ as in
$\oplus_1$.  So by the definition of
$I^{\dg}_\theta(\bar{\cP}^*)$ the set $S^\partial_\theta \backslash
S_2 \in I^{\dg}_\theta(\bar{\cP}^*)$ where $S_2 = \{\zeta \in
S^\partial_\theta$: there is $a \in \cP^*_{< \partial}$ such that
$\otp(a) = \theta,\sup(a) = \zeta$ and $a \subseteq e$ hence $\zeta
\in e\}$.

For each $\zeta \in S$, let $a \in \cP^*_{< \partial}$ witness $\zeta
\in S_2$, as in the proof of clause $(c)(\gamma)$ we get that $\zeta
\in \good''_\theta(\bar{\cP})$.  Clearly this suffices for proving
clauses $(d)(\delta),(\varepsilon)$.
\end{PROOF}

\begin{claim}
\label{b8}
Let $\sigma > \partial > \theta$.

\noindent
1) $S^\sigma_\partial \notin \check I^{\ac}_\partial \langle
   \sigma,\sigma \rangle$ moreover $\check I^{\ac}_\partial \langle
   \sigma,\sigma\rangle$ is a normal ideal on $S^\sigma_\partial$.

\noindent
2) If $S_1 \in \check I_\theta[\sigma]$ and $S_2 \in 
\check I^{\ac}_\theta \langle \sigma,\partial \rangle$ \then \, $S_1 \backslash
   S_2$ is non-stationary.
\end{claim}

\begin{remark}
If $\sigma = \partial^+$, see \ref{m19}.
\end{remark}

\begin{PROOF}{\ref{b8}}
1) Easy.

\noindent
2) Let $\bar{\cP}' = \langle \cP'_\varepsilon:\varepsilon <
\sigma\rangle$ be a $(\sigma,\partial,< \sigma)$-system witnessing
$S_1 \in \check I_\theta[\sigma]$.

Now instead of choosing $N_\varepsilon$ for $\varepsilon \le \partial$
we choose $N_\varepsilon$ and $\bar N_\varepsilon$ 
by induction on $i < \sigma$ such that:
\mn
\begin{enumerate}
\item[$\oplus(A)$]  $(a) \quad N_\varepsilon \prec
  (\cH(\chi),\in,<^*_\chi)$
\sn
\item[${{}}$]  $(b) \quad \|N_\varepsilon\| < \sigma$ and
  $N_\varepsilon \cap \sigma \in \sigma$
\sn
\item[${{}}$]  $(c) \quad \langle N_\zeta:\zeta \le \xi\rangle \in
  N_\varepsilon$ for $\xi < \varepsilon$
\sn
\item[$(B)$]  $(a) \quad \bar N_\varepsilon  = \langle N_{\varepsilon,a}:a \in
  \cP'_\varepsilon \rangle$
\sn
\item[${{}}$]  $(b) \quad N_{\varepsilon,a} \prec (\cH(\chi),\in,<^*_\chi)$
\sn
\item[${{}}$]  $(c) \quad \|N_{\varepsilon,a}\| < \partial$
\sn
\item[${{}}$]  $(d) \quad$ if $a \in \cP'_\varepsilon$ then
$\langle N_{\xi,a \cap \xi}:\xi \in a \cup \{\varepsilon\}\rangle$ is
  $\prec$-increasing and 

\hskip25pt  $\xi \in a \cup \{\zeta\} \wedge \xi = \sup(a \cap \xi) 
\Rightarrow N_{\xi,a \cap \xi} =
  \cup\{N_{\zeta,a \cap \zeta}:\zeta \in a\}$
\sn
\item[${{}}$]  $(e) \quad E,\bar M,\bar\delta,\sigma,\bar{\cP}^*$ and
  $\bar{\cP}'$ belongs to $N_{i,a}$
\sn
\item[${{}}$]  $(f) \quad \langle N_{\zeta,b}:\zeta \le \xi,b \in
  \cP'_\zeta\rangle$ and $\langle N_\zeta:\zeta \le \xi\rangle$
  belongs to $N_{\varepsilon,a}$ and to $N_\varepsilon$

\hskip25pt  when $\xi < \varepsilon_* < \sigma$
\sn
\item[${{}}$]  $(g) \quad \partial +1 \subseteq N_{\varepsilon,a}$.
\end{enumerate}
\mn
The rest should be clear.
\end{PROOF}

\begin{PROOF}{\ref{y1}}
\underline{Proof of \ref{y1}}
1) As $\partial,\theta$ are regular cardinals and $\partial >
\theta^+$ let $\bar{\cP}^* := \langle
\cP^*_\alpha:\alpha < \partial\rangle$ be a $(\partial,\le \theta,<
\partial)$-system satisfying $S^\sigma_\theta \notin
I^{\cg}_\theta(\bar{\cP}^*)$, see \ref{m14}(1), \ref{m19}(3).  
Let $\chi,M_*$ be as in \ref{b4} for our $\lambda$ such
that $\bar{\cP}^* \in M_*$.  Let
$E,\bar{\cP}$ be as constructed in \ref{b4} for our
$\lambda,M_*$ and recall $\alpha \in \nacc(E) \Rightarrow 
\cf(\alpha) = \aleph_0$.  So if $\delta \in E
\cap S^\lambda_\sigma$ \then \, $\delta \in \acc(E)$ and so there is
an increasing continuous sequence $\langle \delta_i:i < \sigma\rangle$
of members of $E$ with limit $\delta$; hence by clauses
$(c)^+(\gamma)$ we have $(\exists^{\stat} i < \delta)
[i \in \good''_\theta(\bar{\cP})]$.

As we have started with any $\delta \in E \cap S^\lambda_\theta$
clearly $\good''_\theta(\bar{\cP})$ reflects in any $\delta \in E \cap
S^\lambda_\sigma$, but $\good''_\theta(\bar{\cP}) \in \check
I_\theta[\lambda]$.  Now by $\boxplus(b)$ of \ref{b4}
$\delta \in S^\lambda_\theta \backslash E \Rightarrow
\delta \in \good''_\theta(\bar{\cP})$ so $\good''_\theta(\bar{\cP})
\in \check I_\theta[\lambda]$ is as required.

\noindent
2) Same proof.

\noindent
3) Similarly using clause $(d)(\varepsilon)$ of \ref{b4}.
\end{PROOF}

\begin{PROOF}{\ref{y12}}

\noindent
\underline{Proof of \ref{y12}}:

\noindent
1) Let $\chi,\lambda,M_*$ be as the assumption of \ref{b4} such that
in addition $2^{\theta^{+n}} < M_* \cap \lambda$ for every $n$.  Let
$E$ and $\bar{\cP} = \langle \cP_\alpha:\alpha < \lambda\rangle$ be as in the
conclusion of \ref{b4}.

Recalling Definition \ref{m13}(2A), let 
$S_* = \good''_\theta(\bar{\cP}) \subseteq S^\lambda_\theta$, 
so obviously $S_* \in \check I_\theta[\lambda]$ and 
for every $n$ let $S_n = \{\delta:\cf(\delta)
=\theta^{+n}$ and $n=0 \Rightarrow \delta \notin S_*$ and $[n \ge 1
  \Rightarrow \delta \cap 
S^\lambda_\theta \backslash S_*$ is a stationary subset of $\delta]\}$.

Note that by the assumption of part of the theorem
\mn
\begin{enumerate}
\item[$\boxplus_1$]  $S_0$ is a stationary subset of $\lambda$.
\end{enumerate}
\mn
For $n \ge 1$ and $\delta \in S_n$ we choose $\langle
\gamma_{\delta,\varepsilon}:\varepsilon < \cf(\delta)\rangle$, an
increasing continuous sequence with limit $\delta$ and let $s_\delta =
\{\varepsilon < \cf(\delta):\cf(\varepsilon) = \theta$ and
$\gamma_{\delta,\varepsilon} \notin S_*\}$, 
so as $\delta \in S_n$ necessary $s_\delta$ is a stationary
subset of $\theta^{+n}$.

For every stationary $s \subseteq S^{\theta^{+n}}_\theta$ let $S_{n,s}
= \{\delta \in S_n:s_\delta = s\}$, the sequence
$\langle S_{n,s}:s \subseteq S^{\theta^{+n}}_\theta$ is 
stationary$\rangle$ is a
partition of $S_n$ and for some club $E_{n,s} \subseteq E$ of $\lambda$ we
have [$S_{n,s} \cap E_{n,s} = \emptyset \Leftrightarrow
S_{n,s}$ is not stationary] for every such $s$.

Let $E_* = \cap\{E_{n,s}:n \ge 1$ and $s \subseteq \theta^{+n}$ is
stationary$\}$, so as we are assuming $2^{\theta^{+n}} < \lambda$,
clearly $E_*$ is a club of $\lambda$. 

Clearly if $``n \ge 1 \wedge (s \subseteq S^{\theta^{+n}}_\theta$
stationary$) \Rightarrow S_{n,s} \subseteq \lambda$ is not stationary"
then $n=0,S=S_0$ satisfy the desired conclusion.  So assume that $n
\ge 1$ and $s \subseteq \theta^{+n}$ is stationary and $S_{n,s}$ is
stationary.  If $S_{n,s}$ reflects in no $S^\lambda_{\theta^{+m}},m >
n$ we are done, and also if $\refl(S_{n,s}) \cap
S^\lambda_{\theta^{+n+1}}$ reflect in no $S^\lambda_{\theta^{+n}},m >
n+1$, we are done. 

Hence it suffices to prove
\mn
\begin{enumerate}
\item[$\boxplus_2$]  if $n \ge 1,s \subseteq S^{\theta^{+n}}_\theta$ is
stationary and $S_{n,s} \subseteq \lambda$ is stationary,  
$n \ge 2,m \ge n+2$ \then \, $S_{n,s}$ does not reflect
in any $\delta_* \in S^\lambda_{\theta^{+m}} \cap \acc(E_*)$.
\end{enumerate}
\mn
Toward this let $\sigma = \theta^{+m}$ and 
$\bar\delta = \langle \delta_i:i < \sigma\rangle$
be an increasing continuous sequence of ordinals from $E_*$ with 
limit $\delta_{i(\sigma)} := \delta_*$.  As $s
\subseteq S^{\theta^{+n}}_\theta$ is stationary and $n \ge 2$, let
$\partial = \theta^{+n}$ by \ref{m6}, \ref{m19}(3) 
there is $\bar{\cP}^* = \langle \cP^*_\zeta:\zeta <
\partial \rangle$ a $(\partial,\theta)$-system such that
$S^\lambda_\theta \notin I^{\cg}_\theta(\bar{\cP}^*)$.

Note that $\bar{\cP}^* \in M_*$ because 
$2^{\theta^{+n}} < \lambda$ and $M_* \cap \lambda$. 
So our $\bar{\cP}$ satisfies the conclusion of \ref{b4}, 
so $\boxplus$ holds indeed hence we are done.

\noindent
2),3) The proof is really included in the proof of part (1).
\end{PROOF}

\begin{remark}
%\label{}
In the proof of \ref{b4}, for regular
$\kappa \in (\theta,\lambda)$ and $s$ a stationary subset of
$S^\kappa_\theta$ we can let $S_{\kappa,s} = \{\delta \in
S^\lambda_\kappa$: for some increasing continuous sequence $\langle
\alpha_i:i < \kappa\rangle$ of ordinals with limit $\delta$, the set
$\{i \in S^\kappa_\theta:i \in s$ iff $\alpha_i \in S_*\}$ is not
stationary$\}$.  Let $E_{\kappa,s}$ be a club of $\lambda$, disjoint
to $S_{\kappa,s}$ if $S_{\kappa,s}$ is not stationary.  Let $\kappa_*
< \lambda$ and $E_* = \cap\{E_{\kappa,s}:\kappa \in (\theta,\kappa_*)$
is regular and $s \subseteq \kappa\}$.  We can then continue as above.
\end{remark}
\bigskip

\subsection {Quite free witnesses of $\pcf$-cases exist} \label{Quite} \
\bigskip

\begin{definition}
\label{b18}
1) We say $(\lambda,\bar\lambda,J,\kappa)$ is a $\pcf$-case (may omit
$J$ when $J = [\kappa]^{< \kappa}$ \when \,:
\mn
\begin{enumerate}
\item[$(a)$]  $\bar\lambda = \langle \lambda_i:i <
  \kappa\rangle$ is a sequence of regular cardinals $> \kappa$
\sn
\item[$(b)$]  $J$ is an ideal on $\kappa$ 
\sn
\item[$(c)$]  $\lambda = \tcf(\prod\limits_{i < \kappa}
  \lambda_i,<_J)$.
\end{enumerate}
\mn
2) We say $\bar f$ witness a $\pcf$-case
$(\lambda,\bar\lambda,J,\kappa)$ or is a witness for it 
\when \, $\bar f$ is $<_J$-increasing and
$<_J$-cofinal in $(\prod\limits_{i < \kappa} \lambda_i,<_J)$.

\noindent
3) We say $\bar f$ obeys $(\lambda,\bar\lambda,J,\bar{\cP},\kappa)$ 
\when \, for some $\bar g,\bar f$ obeys $(\lambda,\bar\lambda,J,
\kappa,\bar{\cP})$ as witnessed by $\bar g$, 
see part (4) below and $\bar f$ witnesses the pcf-case
$(\lambda,\bar\lambda,J,\kappa)$.  Not mentioning $\bar g$ means for
   some $\bar g$.

\noindent
4) We say that $\bar f$ obeys $(\lambda,\bar\mu,J,\kappa,\bar{\cP})$
   as witnessed by $\bar g$ \when \, :
\mn
\begin{enumerate}
\item[$(a)$]  $\bar f = \langle f_\alpha:\alpha < \lambda\rangle$; 
\sn
\item[$(b)$]  $J$ is an ideal on $\kappa$ and $\bar\mu = \langle \mu_i:i <
  \kappa\rangle$
\sn
\item[$(c)$]  $f_\alpha \in {}^\kappa\Ord$
\sn
\item[$(d)$]  $\bar f$ is $<_J$-increasing
\sn
\item[$(e)$]  $\bar{\cP} = \langle \cP_\alpha:\alpha < \lambda\rangle$
  is a $(\lambda,\lambda,\le 2^\lambda)$-system, so \wilog \,
  $\subseteq$-increasing 
\sn
\item[$(f)$]  $\bar g = \langle g_a:a \in \bigcup\limits_{\alpha}
  \cP_\alpha\rangle$
\sn
\item[$(g)$]  $g_a \in {}^\kappa \Ord$
\sn
\item[$(h)$]  $g_a(i) < g_b(i)$ \when \, $a \triangleleft b$ are from
$\cP_{< \lambda}$ and $|b| < \mu_i$ where $\cP_{< \alpha} :=
  \cup\{\cP_\beta:\beta < \alpha\}$
\sn
\item[$(i)$]  if $a \in \cP_\alpha$ then $g_a <_J f_\alpha$
\sn
\item[$(j)$]  if $\beta \in a \in 
\cP_\alpha,i < \kappa$ and $|a| < \mu_i$ then $f_\beta(i) < g_a(i)$.
\end{enumerate}
\end{definition}

\begin{convention}
\label{b20}
We may allow $\bar f = \langle f_\alpha:\alpha \in S\rangle$ where $S
\subseteq \lambda = \sup(S)$, that is, say $\bar f$ obeys
$(\lambda,\bar\mu,J,\kappa,\bar{\cP})$ as witnessed by some $\bar g$ when
$\langle f'_\alpha:\alpha < \lambda\rangle$ satisfies the demands
there where
$\alpha \in S \Rightarrow f'_{\otp(S \cap \alpha)} = f_\alpha$.
\end{convention}

\begin{claim}
\label{b21}
Assume $(\lambda,\bar\lambda,J,\kappa)$ is a pcf-case, $\mu =
\lim\inf_J(\bar\lambda)$ and $\bar{\cP}$ is a 
$(\lambda,\mu,< \lambda)$-system.

\noindent
1) There is $\bar f$ obeying $(\lambda,\bar\lambda,J,\kappa,\bar{\cP})$.

\noindent
2) For every $\bar f$ witnessing $(\lambda,\bar\lambda,J,\kappa)$, for
   some unbounded $S \subseteq \lambda,\bar f \rest S$ obeys
$(\lambda,\bar\lambda,J,\kappa,\bar{\cP})$. 

\noindent
3) If $\bar f$ obeys $(\lambda,\bar\lambda,J,\kappa,\bar{\cP})$ and
   $\theta = \cf(\theta) < \lim\inf_J(\bar\lambda)$ \then \,
$S_{\ged}[\bar f] \supseteq \good''_\theta(\bar{\cP})$.
\end{claim}

\begin{remark}
The proof is like the ones in \cite[Ch.I]{Sh:e}, \cite{Sh:509}.
\end{remark}

\begin{PROOF}{\ref{b21}}
1) Follows by (2).

\noindent
2) Let $\bar f = \langle f_\alpha:\alpha < \lambda\rangle$ witness
the pcf-case $(\lambda,\bar\lambda,J,\kappa)$.

By induction on $\beta < \lambda$ we choose $\langle g_a:a \in
\cP_\beta\rangle$ and $\alpha(\beta)$ such that
\mn
\begin{enumerate}
\item[$\boxplus$]  $(a) \quad g_a \in \Pi \bar\lambda$
\sn
\item[${{}}$]  $(b) \quad$ if $i < \kappa,b \triangleleft a$ and
  $\{a,b\} \subseteq \cP_{< \beta}$
  and $|a| < \lambda_i$ then $g_b(i) < g_a(i)$
\sn
\item[${{}}$]  $(c) \quad \alpha(\beta) < \lambda$ and $\beta_1 <
  \beta \Rightarrow \alpha(\beta_1) < \alpha(\beta)$
\sn
\item[${{}}$]  $(d) \quad$ if $i < \kappa,\beta_1 \in a \in \cP_\beta$
  and $|a| < \lambda_i$ \then \, $f_{\alpha(\beta_1)}(i) < g_a(i)$
\sn
\item[${{}}$]  $(e) \quad$ if $a \in \cP_{\le\beta}$ \then \, $g_a <_J 
f_{\alpha(\beta)}$.
\end{enumerate}
\mn
In stage $\beta$ we first choose $g_a$ for $a \in \cP_\beta
\backslash \cP_{< \beta}$, note that this means that for every $i < \kappa$,
we have to choose $g_a(i)$ as an ordinal $< \lambda_i$, which  is a
regular cardinal and if $|a| < \lambda_i$ it should be bigger than
$\le |a|$ ordinals $< \lambda_i$, so this is easy.

As for $\alpha(\beta)$ for each $a \in \cP_{\le\beta}$, as $\bar f$ is
cofinal in $(\Pi \bar\lambda,<_J)$ there is $\gamma_{\bar a} <
\lambda$ such that $g_a <_J f_{\gamma_a}$.  So $\alpha(\beta)$ should
be an ordinal $< \lambda$ and $> \sup\{\alpha(\beta_1);\beta_1 <
\beta\}$ which is an ordinal $< \lambda$, as $\lambda$ is regular and
it also should be $> \sup\{\gamma_a:a \in \cP_{\le\beta}\}$ which is $<
\lambda$ as $\lambda$ is regular $> |\cP_\alpha|$.

\noindent
3) Straight.
\end{PROOF}

\begin{definition}
\label{b24}
Let $J$ be an ideal on $\kappa$, we may omit it below when $J =
J^{\bd}_\kappa$. 

\noindent
1) A set $\cF \subseteq {}^\kappa\Ord$ is $J$-free \when \, there is a
   sequence $\langle a_f:f \in \cF\rangle$ of members of $J$ such that
   $f_1 \ne f_2 \wedge \{f_1,f_2\} \subseteq \cF \wedge i \in \kappa
   \backslash a_{f_1} \backslash a_{f_2} \Rightarrow f_1(i) \ne f_2(i)$.

\noindent
2) A set $\cF \subseteq {}^\kappa\Ord$ is $(\theta,J)$-free \when \,
$\cF'$ is $J$-free whenever $\cF' \subseteq \cF$ has cardinality
   $< \theta$.

\noindent
3) A sequence $\langle f_\alpha:\alpha < \alpha_*\rangle$ of members
   of ${}^\kappa\Ord$ is a $(\theta,J)$-free sequence \when , for every $u \in
   [\alpha_*]^{< \theta}$ there is a sequence $\langle a_\alpha:\alpha
   \in u\rangle$ of members of $J$ such that: if $\alpha < \beta$ are
   from $u$ then $i \in \kappa \backslash a_\alpha \backslash a_\beta
   \Rightarrow f_\alpha(i) < f_\beta(i)$.

\noindent
4) A set $\cF \subseteq {}^\kappa\Ord$ (we may use a sequence listing
   it) is called $(\theta_2,\theta_1,J)$-free \when \, for every $\cF'
   \subseteq \cF$ of cardinality $< \theta_2$, we can find a partition
   $\langle \cF'_\varepsilon:\varepsilon < \varepsilon(*)\rangle$ of
   $\cF'$ such that:
\mn
\begin{enumerate}
\item[$\bullet$]  each $\cF'_\varepsilon$ has cardinality $< \theta_1$
\sn
\item[$\bullet$]  we can find a sequence $\langle s_f:f \in
  \cF'\rangle$ of members of $J$ such that $f_1 \in
  \cF'_{\varepsilon_1} \wedge f_2 \in \cF'_{\varepsilon_2} \wedge
  \varepsilon_1 \ne \varepsilon_2 \wedge i \in \kappa \backslash
  s_{f_1} \backslash s_{f_2} \Rightarrow f_1(i) \ne f_2(i)$.
\end{enumerate}
\mn
4A) A set $\cF \subseteq {}^\kappa\Ord$ is called $\langle
\theta_2,\theta_1,J \rangle$-free when for every $\cF' \subseteq
   \cF$ of cardinality $\theta_2$, there is a $J$-free $\cF''
   \subseteq \cF'$ of cardinality $\theta_1$.

\noindent
4B) Similarly to 4), 4A) for a sequence $\langle f_\alpha:\alpha <
\alpha_*\rangle$ of members of ${}^\kappa\Ord$ means that it is with
no repetitions and $\{f_\alpha:\alpha \in u\}$ satisfies the requirement.

\noindent
5) A set $\cF \subseteq {}^\kappa\Ord$ is called $\langle
   \theta_2,\theta_1,J\rangle$-stable \when \, for every $u \subseteq
   \Ord$ of cardinality $< \theta_1$ the set $\{f \in \cF:i$ the set
   $\{i < \kappa:f(i) \in u\}$ is not in $J\}$ has cardinality $<
   \theta_2$.

\noindent
5A) A set $\cF \subseteq {}^\kappa\Ord$ is $(\theta,J)$-stable when it
is $(\theta,\theta,J)$-stable.

\noindent
5B) A set $\cF \subseteq {}^\kappa \Ord$ 
is $(\theta_2,\theta_1,J)$-stable \when \, for every
$\theta \in [\theta_2,\theta_1)$ is $(\theta,J)$-stable.
\end{definition}

\noindent
Toward proving Theorem \ref{y19} we prove
\begin{claim}
\label{b40}
If (A) then (B) where:
\mn
\begin{enumerate}
\item[$(A)$]  $(a) \quad (\lambda,\bar\lambda,J,\kappa)$ is a pcf-case
\sn
\item[${{}}$]  $(b) \quad M_* \prec
  (\cH(\lambda^+),\theta,<^*_{\lambda^+})$ has cardinality $<
  \lambda,M_* \cap \lambda \in \lambda$ and
  $(\lambda,\bar\lambda,J,\kappa) \in$

\hskip25pt $M_*$; (clearly exists and by
\ref{b4}, \ref{b21} there are $\bar{\cP},E,\bar f$, as required 

\hskip25pt below) 
\sn
\item[${{}}$]  $(c) \quad \bar f,\bar{\cP},E$ are as in \ref{b4} for our
 $\lambda,M_*$ 
\sn
\item[${{}}$]  $(c) \quad \bar f^1$ obeys
  $(\lambda,\bar\lambda,J,\kappa,\bar{\cP})$
\sn
\item[${{}}$]  $(d) \quad \mu$ is a limit uncountable cardinal
\sn
\item[${{}}$]  $(e) \quad \mu = \lim \inf_J(\bar\lambda)$, i.e. $\mu =
  \min\{\chi$: the set $\{i < \kappa:\lambda_i < \chi\}$ is not 

\hskip25pt from $J\}$
\sn
\item[${{}}$]  $(f) \quad \partial = \cf(\partial) <
  \kappa,J$ is $\partial^+$-complete 
\sn
\item[${{}}$]  $(g) \quad S \subseteq S^\lambda_\partial$ is stationary
  such that $\delta \in S \Rightarrow (\mu^2$ divide $\delta$)
\sn
\item[${{}}$]  $(h) \quad \bar \alpha = \langle
  \alpha_{\delta,i}:\delta \in S,i < \partial \rangle$ where
  $\bar\alpha_\delta = \langle \alpha_{\delta,i}:i < \partial \rangle$ is
increasing

\hskip25pt  continuous with limit $\delta$ such that
  $\alpha_{\delta,i}$ is divisible by $\mu$
\sn
\item[${{}}$]  $(i) \quad \bar f = \bar f^2 = \langle f^2_\delta:\delta \in
  S\rangle$ is where $f^2_\delta:\partial \times \kappa \rightarrow
  \delta$ is defined by $f^2_\delta(i,j) =$

\hskip25pt $\alpha_{\delta,i} + f^1_\delta(j)$
\sn
\item[${{}}$]  $(j) \quad J_* = J^{\bd}_\partial \times J = \{u
 \subseteq \partial \times \kappa$: for every $i < \partial$ large enough,
$\{j < \kappa:$

\hskip25pt $(i,j) \in u\} \in J\}$; of course, we can translate
  $J_*$ to an ideal 

\hskip25pt on $\{v \subseteq \kappa:\{(i,j) \in \partial
  \times \kappa: \partial \cdot j + i \in v\} \in J_*\}$.
\sn
\item[$(B)$]  $(a)(\alpha) \quad$ if $\theta \in [\kappa,\mu)$ then
the sequence $\bar f^2$ is 
$(\theta^{+\comp(J)+1},\theta^{+4},J_*)$-free 

\hskip35pt  recalling $\partial < \comp(J) 
\le \kappa$, see \ref{b26} and \ref{m0}(5)
\sn
\item[${{}}$]  \hskip10pt $(\beta) \quad \bar f^2$ is
  $(\comp(J),J_*)$-free
\sn
\item[${{}}$]  \hskip10pt $(\gamma) \quad$ if 
$\theta \in [\kappa,\mu)$ is a limit cardinal
and $\cf(\theta) \notin [\comp(J),\kappa^+)$ and 

\hskip35pt $(\forall \Upsilon)(\kappa < \Upsilon < \mu \wedge \cf(\Upsilon) \in
[\comp(J),\kappa^+) \Rightarrow \pp_J(\mu) < \theta)$ 

\hskip35pt  \then \, $\bar f^2$ is $(\theta^{+\comp(J)+1},\theta^+,J_*)$-free
\sn 
\item[${{}}$]  $(b) \quad$ if $\sigma$ is regular and $\delta \in
  S^\lambda_\sigma$ and $\sigma < \mu$ \then \,, see Definition \ref{m60}:
\sn
\item[${{}}$]  $\quad (\alpha) \quad \kappa^{+4} \le \sigma
  \Rightarrow \delta \notin S^{\ch}_J[\bar f]$
\sn
\item[${{}}$]  $\quad (\beta) \quad \kappa^+ < \sigma <
  \kappa^{+\comp(J)+1} \Rightarrow \delta \notin S^{\bd}_J[\bar f]$
\sn
\item[${{}}$]  $\quad (\gamma) \quad \kappa \le \theta \wedge
  \theta^{+4} \le \sigma < \theta^{+\comp(J)+1} \Rightarrow \delta
  \notin S^{\bd}_J[\bar f]$.
\end{enumerate}
\end{claim}

\begin{remark}
\label{b41}
This continues \cite{Sh:108} and \cite{MgSh:204}; note that 
here $\partial < \kappa$.  This helps; there are relatives with
$\sigma \ge \kappa$ but not needed at present.
\end{remark}

\begin{PROOF}{\ref{b40}}
Note that
\mn
\begin{enumerate}
\item[$\boxplus_1$]  if $\theta = \cf(\theta) \in \mu \backslash
  \kappa^+$ \then \, $S_{\ged}[\bar f] \cap S^\lambda_\theta \supseteq
  \good''_\theta[\bar{\cP}]$.
\end{enumerate}
\mn
[Why?  By \ref{b21}(3).]
\mn
\begin{enumerate}
\item[$\boxplus_2$]  if $\theta,\sigma$ are regular cardinals from
  $(\kappa,\mu)$ and $\theta^{+2} < \sigma$ \then \, 
$S_{\ged}[\bar f] \cap S^\lambda_\theta$ reflect in every $\delta \in
  S^\lambda_\sigma$. 
\end{enumerate}
\mn
[Why?  Let $\Upsilon = \theta^{+2}$, hence by \ref{m19}(3)
 there is a $(\Upsilon,\theta,< \Upsilon)$-system 
such that $S^\Upsilon_\theta \notin I^{\cg}_\theta[\Upsilon]$, 
see Definition \ref{m14}(1) 
hence by \ref{b4}, that is the choice of
  $\bar{\cP}$, the set $\good''_\theta(\bar{\cP}) \subseteq
  S^\lambda_\theta$ reflect in every $\delta \in S^\lambda_\sigma$,
  and so by $\boxplus_1$ we are done.]
\mn
\begin{enumerate}
\item[$\boxplus_3$]  $S^{\ged}_J[\bar f]$ include $\{\delta <
  \lambda:\theta^{+4} \le \cf(\delta) < \theta^{+\comp(J)+1}\}$ when
$\theta \in [\kappa,\mu)$.
\end{enumerate}
\mn
[Why?  By $\boxplus_2$, \ref{m62}(2), \ref{m65}(1),(3).]

So we have proved (b) of (B); concerning $(B)(b)(\gamma)$ recall that
\mn
\begin{enumerate}
\item[$\bullet$]  if $\delta \in S^{\ch}_I[\bar f]$ then for some club
  $e$ of $\delta$ we have $\alpha \in e \wedge \cf(\alpha) > \kappa
  \Rightarrow \alpha \in S^{\ch}_J[\bar f]$, (similarly for
  $S^{\ged}_J[\bar f]$)
\sn
\item[$\boxplus_4$]  $\bar f^2$ is
  $(\kappa^{+\comp(J)+1},\kappa^{+4},J)$-free, see Definition
\ref{b24}(4), that is as a set.
\end{enumerate}
\mn
[Why?  By $\boxplus_6$ proved below using $\boxplus_3$.]
\mn
\begin{enumerate}
\item[$\boxplus_5$]  if $\theta \in [\kappa,\mu)$ then $\bar f^2$ is 
$(\theta^{+\comp(J)+1},\theta^{+4},J)$-free.
\end{enumerate}
\mn
[Why?  By $\boxplus_6$ below using $\boxplus_3$.]
\mn
\begin{enumerate}
\item[$\boxplus_6$]  if $\theta_2 > \theta_1 = \cf(\theta_1) > \kappa$
and $\delta < \lambda \wedge \theta_1 \le \cf(\delta) < \theta_2
\Rightarrow \delta \in S^{\ged}_J[\bar f]$ \then \, $\bar f^2$ 
is $(\theta_2,\theta_1,J_*)$-free.
\end{enumerate}
\mn
Toward this we prove for $\theta \in [\theta_1,\theta_2)$ that
\mn
\begin{enumerate}
\item[$\oplus_{\bar f,\theta}$]  if $u \subseteq S$, recalling $S
  \subseteq S^\lambda_\sigma,|u|
= \theta$ \then \, we can find $\bar s = \langle s_\alpha:\alpha \in u
\rangle \in {}^u (J_*)$ such that in the graph $(u,R_{\bar s})$ every node
has valency $< \theta_1$ where:
\sn
\begin{enumerate}
\item[$\bullet$]  for $u \subseteq \lambda$ and $\bar s \in {}^u J_*$
  let $(u,R_{\bar s})$ be the following graph: $\alpha R_{\bar s}
  \beta$ iff $\alpha \ne \beta \in u$ and for some $(i,j) \in \sigma
  \times \kappa$, we have $(i,j) \notin s_\alpha \cup s_\beta$ and
  $f_\alpha(i) = f_\beta(i)$.
\end{enumerate}
\end{enumerate}
\mn
\underline{Why this suffice}?  As then let $\langle u_t:t \in
I \rangle$ list the components of the graph $(u,R_{\bar s})$, so
necessarily each component has cardinality $< \theta$, recalling
$\theta_1$ is regular, so $\langle \{f_\alpha:\alpha \in u_t\}:t \in
I \rangle$ is a partition as required in Definition \ref{b24}(4).

\noindent
\underline{Why this is true}?  We prove this by induction on
$\otp(u)$.
\medskip

\noindent
\underline{Case 1}:  $\otp(u) < \theta_1$

Let $s_\alpha = \emptyset \in J_*$ for $\alpha \in u$, clearly as required.
\medskip

\noindent
\underline{Case 2}:  $\otp(u) = \zeta +1$

Let $\alpha = \max(u)$, let $\bar s \in {}^{u \cap \alpha} (J_*)$ be as
promised for $u \cap \alpha$ and define $\bar s' \in {}^u (J_*)$ by
$s'_\beta$ is $s_\beta$ if $\beta < \alpha$ and is $\emptyset$ if
$\beta = \alpha$, now check.
\medskip

\noindent
\underline{Case 3}:  $\delta = \otp(u)$ is a limit ordinal of
cofinality $< \theta_1$

Let $\sigma := \cf(\delta)$ and $\langle
\alpha_\varepsilon:\varepsilon < \sigma\rangle$ be increasing
continuous with limit $\sup(u)$ such that $\alpha_0=0$.  For $\varepsilon <
\sigma$ let $u_\varepsilon = u \cap [\alpha_\varepsilon,\alpha_{\varepsilon
+1})$ and let $\bar s_\varepsilon = \langle s_\alpha:\alpha \in 
u_\varepsilon \rangle$ be as required for 
$u_\varepsilon$, exists as $\otp(u_\varepsilon) < \otp(u)$.  
So $\bar s = \langle s_\alpha:\alpha \in u\rangle$ is well defined.
Now for each $\beta \in u,(i_*,j_*) \in \partial \times \kappa$ and
$\varepsilon$ the set $w_{\beta,\varepsilon,i_*,j_*} = \{\gamma \in
u_\varepsilon:(i_*,j_*) \notin s_\gamma$ and $f^2_\gamma(i_*,j_*) =
f^2_\beta(i_*,j_*)\}$ has cardinality $< \theta_1$ because
$\gamma_1,\gamma_2 \in w_{\beta,\varepsilon,i_*,j_*} \Rightarrow
((i_*,j_*) \in \partial \times \kappa \backslash (s_{\gamma_1} \cup
s_{\gamma_2})) \wedge f^2_{\gamma_1}(i_*,j_*) =
f^2_{\gamma_2}(i_*,j_*)$;  hence $w_\beta :=
\cup\{w_{\beta,\varepsilon,i,j}:\varepsilon < \otp(u)$ and $i <
\partial,j < \kappa\}$ has cardinality $< \theta_1$ and $\langle
w_\beta:\beta \in u \rangle$ is as required.
\medskip

\noindent
\underline{Case 4}:  $\delta = \otp(u)$ has cofinality $\ge \theta_1$.

We choose $\bar s \in {}^u (J_*),\bar\beta,\bar a^1$ such that:
\mn
\begin{enumerate}
\item[$(*)_{6.1}$]  $(a) \quad \bar\beta = \langle
\beta_\varepsilon:\varepsilon < \cf(\delta)\rangle$ is increasing
continuous
\sn
\item[${{}}$]  $(b) \quad \beta_0 = 0$
\sn
\item[${{}}$]  $(c) \quad \cup\{\beta_\varepsilon:i < \cf(\delta)\} =
\sup(u)$
\sn
\item[${{}}$]  $(d) \quad \bar a^1 = \langle
a^1_\varepsilon:\varepsilon < \cf(\delta)$ non-limit$\rangle$
\sn
\item[${{}}$]  $(e) \quad a^1_\varepsilon \in J$
\sn
\item[${{}}$]  $(f) \quad$ if $\varepsilon > 0$ then
$\beta_\varepsilon = \sup(u \cap \beta_\varepsilon)$
\sn
\item[${{}}$]  $(g) \quad$ if $\varepsilon,\zeta < \cf(\delta)$ are
non-limit and $j \in \kappa \backslash a^1_\varepsilon \backslash
a^1_\zeta$ then $f^1_{\beta_\varepsilon}(j) < f^1_{\beta_\zeta}(j)$
\sn
\item[${{}}$]  $(h) \quad \beta_\varepsilon \in S^\lambda_\sigma$ iff 
$\cf(\varepsilon) = \sigma$.
\end{enumerate}
\mn
[Why such $\bar\alpha,\bar a$ exist?  First, $\sup(u) 
\in S^{\ged}_J[\bar f^1]$ holds by an assumption of $\boxplus_6$
because $\theta_1 \le \cf(\sup(u))$ by the case assumption and
$\cf(\sup(u)) < \theta_2$ as $|u| \le \theta_2$.  Second, use
Definition \ref{m60}(1) recalling clause (d) of $(*)_{6.1}$.]
\mn
\begin{enumerate}
\item[$(*)_{6.2}$]  we can find $\bar a$ such that:
\sn
\item[${{}}$]  $(a) \quad \bar a = \langle a_\varepsilon:\varepsilon <
\cf(\delta)\rangle$
\sn
\item[${{}}$]  $(b) \quad a_\varepsilon = a^1_\varepsilon$ if
$\varepsilon$ is non-limit
\sn
\item[${{}}$]  $(c) \quad a_\varepsilon \in J$
\sn
\item[${{}}$]  $(d) \quad$ if $\varepsilon < \zeta < \cf(\delta)$ and
$\cf(\zeta) < \comp(J)$ or $\cf(\zeta) > \kappa$ then 

\hskip25pt $j \in \kappa \backslash a_\varepsilon 
\backslash a_\zeta \Rightarrow
f_{\beta_\varepsilon }(j) < f_{\beta_\zeta}(j)$.
\end{enumerate}
\mn
[Why?  For non-limit $\varepsilon < \cf(\delta)$ let $a_\varepsilon =
a^1_\varepsilon $.

If $\varepsilon  < \cf(\delta)$ and $\aleph_0 \le \cf(\varepsilon) <
\comp(J)$ then let $e_\varepsilon$ be an unbounded subset of
$\varepsilon$ of order type $\cf(\varepsilon)$ and let $a_\varepsilon
= \kappa \backslash \{i < \kappa:i \notin \cup\{a_{\beta_{\zeta
+1}}:\zeta \in e_\varepsilon\}$ and $f^1_{\beta_\varepsilon}(i) <
f^1_{\beta_{\varepsilon +1}}(i)$ and $\zeta \in e_\varepsilon
\Rightarrow f^1_{\beta_{\zeta +1}}(i) < f^1_{\beta_\varepsilon}(i)\}$.

As $J$ is $\comp(J)$-complete ideal on $\kappa$ and $\bar f^1$ is
$<_J$-increasing clearly $a_\varepsilon \in J$.

If $\varepsilon < \cf(\delta)$ and $\cf(\varepsilon ) > \kappa$ then
let $a_\varepsilon = \{i < \kappa$: the set $\{\zeta < \varepsilon:i
\notin a_{\zeta +1}$ and $f_{\beta_{\zeta +1}}(i) <
f_{\beta_\varepsilon}(i)\}$ is a bounded subset of $\varepsilon\}$.

Toward proving $a_\varepsilon \in J$, first we find $\xi(\varepsilon)
< \varepsilon$ such that: if $i < \kappa$ and the set $\{\zeta <
\varepsilon:i \in \kappa \backslash a_{\zeta +1}$ and $f_{\beta_{\zeta
+1}}(i) < f_{\beta_\zeta}(i)\}$ is bounded below $\varepsilon$ then it
is $\le \xi(\varepsilon)$; this is possible as $\cf(\varepsilon)
> \kappa$.

So $\kappa \backslash a_\varepsilon \supseteq \{i <
\kappa:f^1_{\beta_{\xi(\varepsilon)+1}} < f^1_{\beta_\varepsilon}(i)$
and $i \notin a_{\xi(\varepsilon)+1}\}$ and the latter set is $=\kappa \mod J$
because $(a_{\xi(\varepsilon)+1} \in J) \wedge
(f_{\beta_{\xi(\varepsilon)+1}} <_J f^1_{\beta_\varepsilon})$; it
follows that $a_\varepsilon \in J$.

In the remaining cases $\cf(\varepsilon) \in [\comp(J),\kappa]$ let
$a_\varepsilon  = \kappa \backslash \{i <
\kappa:f_{\beta_\varepsilon}(i) < f_{\beta_{\varepsilon +1}}(i)$ and
$i \notin a_{\varepsilon +1}\}$.  Actually only the $a_\varepsilon$
for $\varepsilon \in S^{\cf(\delta)}_\partial$ are used later.

Let us check that $\langle a_\varepsilon:\varepsilon <
\cf(\delta)\rangle$ is as required in $(*)_{6.2}$ 
so assume $\varepsilon < \zeta <
\cf(\delta)$ and $i \in \kappa \backslash a_\varepsilon \backslash
a_\zeta$.  First, if $\varepsilon,\zeta$ are non-limit then $i \in
\kappa \backslash a^1_\varepsilon \backslash a^1_\zeta$ hence
$f_{\beta_\varepsilon}(i) < f_{\beta_\zeta}(i)$.  Second, if $\varepsilon$ is
non-limit and $\cf(\zeta) < \comp(J)$ then we can find $\xi \in
e_\zeta$ which is $> \varepsilon$, so $i \notin a_{\beta_{\xi +1}}$ as
$a_{\beta_{\xi +1}} \subseteq a_{\beta_\varepsilon}$ hence
$f_{\beta_\varepsilon}(i) < f_{\beta_{\xi +1}}(i)$ and by the choice
of $a_{\alpha_\varepsilon}$ also $f_{\beta_{\xi +1}}(i) <
f_{\beta_\zeta}(i)$, together $f_{\beta_\varepsilon}(i) <
f_{\beta_\zeta}(i)$.  Third, if $\varepsilon$ is a limit ordinal and
$\cf(\zeta) < \comp(J)$ so by the choice of $a_\varepsilon$ we have
$(f_{\beta_\varepsilon}(i) < f_{\beta_{\varepsilon + 1}}(i)) \wedge
(a_\varepsilon \supseteq a_{\varepsilon +1})$ so $i \notin
a_{\varepsilon +1}$; so by the above applied to $(\varepsilon
+1,\zeta)$ we have $f_{\beta_{\varepsilon +1}}(i) <
f_{\beta_\zeta}(i)$, so together $f_{\beta_\varepsilon}(i) <
f_{\beta_\zeta}(i)$.  The cases when $\cf(\zeta) > \kappa$ is
similar.  So we have proved $(*)_{6.2}$.]

Now for each $\varepsilon < \cf(\delta)$ let $u_\varepsilon = u
\cap[\beta_\varepsilon,\beta_{\varepsilon +1})$ hence
$\otp(u_\varepsilon) < \otp(u) = \delta$ hence there is a sequence
$\langle s^\varepsilon_\alpha:\alpha \in u_\varepsilon\rangle$ of members of
$J_*$ as required.  For each $\varepsilon <\cf(\delta)$ and $\beta \in
u_\varepsilon \backslash \{\beta_\varepsilon\}$ let $\bold i(\beta) <
\partial$ be such that $\{\alpha_{\beta,i}:i \in [\bold i(\beta),\sigma)\}
\cap \beta_\varepsilon = \emptyset$ and if $\varepsilon <
\cf(\delta),\beta = \beta_\varepsilon$ so $\beta_\varepsilon \in
S^\lambda_\sigma$ let $\bold i(\alpha) = 0$.

Lastly, let us define $\bar s = \langle s_\beta:\beta \in u\rangle$:
\mn
\begin{enumerate}
\item[$(*)$]  if $\beta \in u_\varepsilon$ then $s_\beta :=
s^\varepsilon_\beta \cup \{(i,j) \in \partial \times \kappa:i \le
\bold i(\beta)\} \cup \{(i,j) \in \partial \times \kappa:j \in
a_\varepsilon \cup a_{\varepsilon +1}\} \cup\{(i,j) \in \partial \times
\kappa:\neg(f^1_{\beta_\varepsilon}(j) \le f^1_\beta(j) <
f^1_{\beta_{\varepsilon  +1}}(j)\}$.
\end{enumerate}
\mn
Let $\beta \in u$ and let $w_\beta = \{\gamma \in u$: 
there is $(i,j) \in \sigma \times \kappa \backslash s_\beta
\backslash s_\gamma$ satisfying $f_\gamma(i,j) = f_\beta(i,j)\}$ and we have
to prove that $w$ has cardinality $< \theta_1$.  Let $\varepsilon <
\cf(\delta)$ be such that $\beta \in u_\varepsilon$ that is $\beta \in
[\beta_\varepsilon,\beta_{\varepsilon +1})$, clearly $\varepsilon$
exists and is unique.  As $s_\beta \supseteq s^\varepsilon_\beta$
clearly $w_\beta \cap [\beta_\varepsilon,\beta_{\varepsilon +1})$ have
cardinality $< \theta_1$.  Now if $\gamma \in u \cap
\beta_\varepsilon \wedge  \beta > \beta_\varepsilon$ then by the
choice of $s_\beta$ we have $s_\beta \supseteq \bold i(\beta) \times
\kappa$ and by the choice of $\bold i(\beta)$ we have 
$\gamma \notin w_\beta$ recalling $\{\alpha_{\gamma,j}:j < \partial\}
\subseteq \beta_\varepsilon$.  If $\gamma \in u \cap \beta_\varepsilon
\wedge \beta = \beta_\varepsilon$ then necessarily
$\beta_\varepsilon \in S^\lambda_\partial$ so $\cf(\beta_\varepsilon) =
\partial$ and let $\xi < \cf(\delta)$ be such that $\gamma \in
[\beta_\xi,\beta_{\xi +1})$, now if $(i,j) \in \partial \times \kappa
\backslash s_\beta \backslash s_\gamma$ \then \, by $(*)_{6.2}(d)$ we have
$f^1_\gamma(i) < f^1_{\alpha_{\xi +1}}(i) < f^1_{\alpha_\varepsilon}(i)$ so
$\gamma \notin w_\varepsilon$.  Together $w_\beta \cap
\alpha_\varepsilon = \emptyset$.

Next, assume $\gamma \in u \backslash \beta_{\varepsilon +1}$ say
$\gamma \in u_\xi,\xi > \varepsilon$; if $\cf(\xi) \ne \partial \vee
\gamma > \beta_\xi$ we use $\bold i(\gamma) \times \kappa \subseteq
s_\gamma$ and if $\cf(\xi) = \partial \wedge \gamma = \beta_\xi$ we use
the chocies of $a_\xi,a_\varepsilon$; hence $w_\beta \backslash
\beta_{\varepsilon +1} = \emptyset$.

Together $w_\beta$ has cardinality $< \theta_1$ as required.  
So we are done proving Case 4, hence proving $\boxplus_6$.
\mn
\begin{enumerate}
\item[$\boxplus_7$]  the sequence $\bar f^2$ is
$(\comp(J)^+,J_*)$-free; this is clause $(a)(\beta)$ of (B).
\end{enumerate}
\mn
[Why?  Let $u \subseteq \lambda$ have cardinality $\le \comp(J)$, let
$\langle \beta_\varepsilon:\varepsilon < |u|\rangle$ list $u$ and
$a_\varepsilon = \{i < \kappa$: for some $\zeta < \varepsilon$ we have
$f^1_{\beta_\zeta}(i) = f^1_{\beta_\varepsilon}(i)\}$, so as $J$ is
$|u|^+$-complete by the assumption clearly $a_\varepsilon \in J$.  Let
$s_{\beta_\varepsilon} = \partial \times a_\varepsilon$ for $\varepsilon
< |u|$,  recalls that for each
$\zeta < \varepsilon,\{i < \kappa:f^1_{\beta_\zeta}(i) =
f^1_{\beta_\varepsilon}(i)\} \in J$ by clause (A)(c) of the assumption
and so $\langle s_\beta:\beta \in u \rangle$ is as required.]
\mn
\begin{enumerate}
\item[$\boxplus_8$]  if $\theta \in [\kappa,\mu)$ then $\bar f^2$ is
  $(\theta^{+\comp(J)+1},\theta^{+4},J_*)$-free.
\end{enumerate}
\mn
[Why?  By $\boxplus_6$ and $(B)(b)(\gamma)$ which we have proved in
  $\boxplus_3$.]
\mn
\begin{enumerate}
\item[$\boxplus_9$]  if $\theta \in [\kappa,\mu)$ is a limit cardinal
  and $\cf(\theta) \notin [\comp(J),\kappa^+)$ and $(\forall
  \Upsilon)(\kappa < \Upsilon < \mu \wedge \cf(\Upsilon) \in
  [\comp(J),\kappa^+) \Rightarrow \pp_J(\mu) < \theta)$ \then \, $\bar f^2$
is $[\theta^{+\comp(J)+1},\theta,J_*)$-free.  This is clause $(B)(a)(\gamma)$
  of the desired conclusion.
\end{enumerate}
\mn
Why?  Clearly $\theta \ne \kappa$ hence recalling $\theta$ is a limit
ordinal $\ge \kappa$ we have $\theta \ge \kappa^{+4}$. 
Again by $\boxplus_6$ it suffices to prove that
if $\delta < \lambda$ and $\cf(\delta) \in
  [\theta,\theta^{+\comp(J)+1})$ then $\delta \notin S^{\ch}_J[\bar f]$
and $\delta \notin S^{\bd}_J[f]$. 

If $\cf(\delta) \ge \theta^{+4}$ this holds by $\boxplus_3$, so we can
assume $\cf(\delta) \in \{\theta^{+ \ell}:\ell \le 3\}$.  Now $\delta
\notin S^{\ch}_J[\bar f]$ as otherwise there is a club $e$ of $\delta$
such that $\alpha \in e \wedge \cf(\alpha) > \kappa \Rightarrow \alpha
\in S^{\ch}_J[f]$, contradicting $\boxplus_3$ applied to
$\kappa^{+4}$.

Also $\delta \notin S^{\bd}_J[\bar f]$ as otherwise $\cf(\delta) =
(\prod\limits_{i < \kappa} \sigma_i,<_J)$ for some $\sigma_i =
\cf(\sigma_i) \in (\kappa,\mu)$ but this contradicts the assumption of
$\boxplus_9$, e.g. $(B)(b)(\gamma)$.
\end{PROOF}

\begin{PROOF}{\ref{y19}}

\noindent
\underline{Proof of \ref{y19}}:

The proof is by cases.
\medskip

\noindent
\underline{Case 1}:  $\lambda$ is singular.

In this case there is a $\mu^+$-free $\cF \subseteq {}^\kappa \mu$ of
cardinality $2^\mu = \lambda$ by \cite[3.10(3)=1f.28(3)]{Sh:898};
more fully by \cite[Ch.II,2.3,pg.53]{Sh:g} for every $\chi \in
(\mu,\lambda)$ there is a $\mu^+$-free $\cF_\chi \subseteq {}^\kappa
\mu$ of cardinality $\chi$; by letting $\bar\chi = \langle
\chi_\varepsilon:\varepsilon < \cf(\lambda)\rangle$ be increaasing
with limit $\lambda$, combining the $\cF_{\chi_\varepsilon}$'s and
$\cF_{\cf(\lambda)}$ we are done.  So clause (A) holds and we are done.
\medskip

\noindent
\underline{Case 2}:  $\lambda$ is regular and $|\alpha|^{< \kappa} =
\lambda$ for some $\alpha < \lambda$.

In this case by \cite[3.6=1f.21]{Sh:898} there is a $\mu^+$-free $\cF \subseteq
{}^\kappa \mu$ of cardinality $2^\mu = \lambda$ so again clause (A)
holds and we are done.
\medskip

\noindent
\underline{Case 3}:  $\lambda$ is regular and $\alpha < \lambda
\Rightarrow |\alpha|^{< \sigma} < \lambda$.

Let $E = \{\delta < \lambda:\alpha < \lambda \Rightarrow |\alpha|^{<
\kappa} < \delta$ and $\delta$ is divisible by $\mu \cdot \mu\}$, clearly a
club of $\lambda$.

Let $S \subseteq E$ be any stationary subset of $S^\lambda_\sigma$.  We choose
$\langle \bar\alpha_\delta:\delta \in S\rangle$ such that
$\bar\alpha_\delta = \langle \alpha_{\delta,i}:i < \sigma\rangle$ is
increasing with limit $\delta$ such that each
$\alpha_{\delta,i}$ is divisible by $\mu$.  By the case assumption we
have $S \in \check I_\sigma[\lambda]$, hence \wilog \,
$\alpha_{\delta_1,i_1} = \alpha_{\delta_2,i_2} \Rightarrow i_1 = i_1
\wedge (\forall i < i_1)(\alpha_{\delta_1,i} = \alpha_{\delta_2,i})$.

Now as $\mu \in \bold C_\kappa$, there is a sequence $\bar\lambda$ such that
$(\lambda,\bar\lambda,J^{\bd}_\kappa,\kappa)$ is a pcf-case such that
$\bar\lambda$ is an increasing sequence of regular cardinals 
with limit $\mu$.  We can choose $\chi,M_*$
as in the assumption of \ref{b4} for $\lambda$ such that $\cH(\mu)
\in M_*$ and the choose $E,\bar{\cP}$ as in the conclusion of \ref{b4}.

Hence by \ref{b21}(1) we can find $\bar f^1 = \langle f^1_\alpha:\alpha <
\lambda\rangle$ obeying
$(\lambda,\bar\lambda,J^{\bd}_\kappa,\kappa,\bar{\cP})$.  Let
$\cd:{}^{\kappa >}\mu \rightarrow \mu$ be one-to-one, we may assume
that $(\forall i)\lambda_i > \kappa$ and $\nu \in
\prod\limits_{j < \kappa} \lambda_j \wedge i < j < \kappa \Rightarrow
\cd(\nu \rest i) < \cd(\nu \rest j)$.  Define
$f^*_\alpha:\kappa \rightarrow \mu$ by $f^*_\alpha(i) = \cd(f_\alpha
\rest (i+1))$, so $f^*_\alpha$ is increasing.

Lastly, let $\alpha_{\delta,i,j} = \alpha_{\delta,i} + f^*_\delta(j)$
and we should prove that $\langle \alpha_{\delta,i,j}:\delta \in S,i <
\sigma,j < \kappa\rangle$ is as required in Definition \ref{y40}, so
$\eta_\delta = \langle \alpha_{\delta,i,j}:(i,j) \in \sigma \times
\kappa\rangle$.  If we have used $f^1_\alpha$ instead of $f^*_\alpha$
we just have to omit clause (d) of \ref{y40}.

Clauses (a),(c) of \ref{y40} holds by our choice of $\eta_\delta$.
Clause (b) of \ref{y40} holds by the choice of $S$ noting that $S \in
\check I_\sigma[\lambda]$ as $S \subseteq E \cap S^\lambda_\sigma$ and
the case assumption.  Clause (d) of
\ref{y40} holds by the choices of the $\bar\alpha_\delta$'s and 
of $\cd,f^*_\alpha$ recalling
$f^1_\alpha \in {}^\kappa \mu$ and $\alpha_{\delta,i}$ is divisible by
$\mu$.  Clause (e) holds by \ref{b40}, that is $(B)(a)$ there says
$\bar f = \bar f^2$ is $(\theta^{+\kappa +1},\theta,J_*)$-free when
$\theta \in [\kappa,\mu)$.  Also clause (f) of \ref{y40} that is
 ``$\bar f$ is $(\kappa^+,J_*)$-free" holds by direct inspection or
 see clause $(B)(a)(\beta)$ of \ref{b40}.

Lastly, clause (g)$'$ follows by clause (g) and clause (g) holds by
\cite{Sh:775}. 
\end{PROOF}

\begin{definition}
\label{b26}
Let $J$ be an ideal on $\kappa$.

\noindent
1) We say $\cF \subseteq {}^\kappa\Ord$ is strongly 
semi-$\langle \theta_2,\theta_1,J \rangle$-stable there are 
no $f_\varepsilon \in \cF$ for $\varepsilon < \theta_2$ 
and $u \subseteq \Ord$ of 
cardinality $< \theta_1$ such that for $\varepsilon < \zeta <
\theta_2$ the following set $A_{\varepsilon,\zeta} =
A_{\kappa,\zeta}(u,\langle f_\varepsilon:\varepsilon \in u\rangle)$
is $\ne \emptyset \mod J$

\[
A_{\varepsilon,\zeta} := 
\{i < \kappa:\min(u \cup \{\infty\} \backslash f_\varepsilon(i)) \ne
\min(u \cup \{\infty\} \backslash f_\zeta(i))\}.
\]

\mn
2) For $<_J$-increasing $\bar f = \langle f_\alpha:\alpha <
\alpha_*\rangle,f_\alpha \in {}^\kappa\Ord$ we say $\bar f$ is 
strongly-semi-$\langle \theta_2,\theta_1,J \rangle$-stable 
(sequence) \when \, there are
no $v \subseteq \alpha_*$ of cardinality $\theta_2$ and $u \subseteq \Ord$ of cardinality
 $< \theta_1$ such that: if $\alpha < \beta$ are from $v$ \then \, the
following set is $\ne \emptyset \mod J$

\[
\{i < \kappa:\min(u \cup \{\infty\} \backslash f_\alpha(i)) \nleq
\min(u \cup \{\infty\} \backslash f_\beta(i)\}.
\]

\mn
3) In parts (1),(2) above, if $\theta_1 = \theta_2$ we may write
$(\theta,J)$ instead of $(\theta_1,\theta_2)$.

\noindent
4) In parts (1),(2) above writing $(\theta_2,\theta_1,J)$ instead of 
$\langle \theta_2,\theta_1,J \rangle$ means: 
strongly-semi-$(\theta,J)$-stable for every $\theta \in [\theta_1,\theta_2)$.
\end{definition}

\begin{claim}
\label{b28}
Assume $\bar f = \langle f_\alpha:\alpha < \lambda\rangle$ witness
the pcf-case $(\lambda,\bar\lambda,J,\kappa)$ and is 
strongly-semi-$(\theta_2,\theta_1,J)$-stable, see \ref{b26}(2),(4).  
\Then \, $S_{\ged}[\bar f] \supseteq
\{\delta < \lambda:\cf(\delta) \in [\theta_1,\theta_2)\}$.
\end{claim}

\begin{PROOF}{\ref{b28}}
Straightforward.
\end{PROOF}

\noindent
Note also
\begin{observation}
\label{b30}
Let $J$ be an ideal on $\kappa$.

\noindent
1) If $f_\alpha \in {}^\kappa\Ord$ for $\alpha < \alpha_*$ and the
   sequence $\langle f_\alpha:\alpha < \alpha_*\rangle$ is
   $(\theta,J)$-free \then \, the set $\{f_\alpha:\alpha < \alpha_*\}$
   is $(\theta,J)$-free and is with no repetitions.

\noindent
2) Similarly for $(\theta_2,\theta_1,J)$-free.

\noindent
2A) Similarly for $\langle \theta_2,\theta_1,J \rangle$-free.

\noindent
3) If $\theta'_2 \ge \theta_2 \ge \theta_1 \ge \theta'_1$ \then \,
\mn
\begin{enumerate}
\item[$(a)$]   $\cF$ is $(\theta_2,J)$-free implies $\cF$ 
is $(\theta_1,J)$-free
\sn
\item[$(b)$]  similarly for $\bar f$
\sn
\item[$(c)$]  $\cF$ is $\langle \theta_2,\theta_1,J\rangle$-stable
implies $\cF$ is $\langle \theta'_2,\theta'_1,J\rangle$-stable.
\end{enumerate}
\mn
4) If $\cF \subseteq {}^\kappa\Ord$ is $(\theta^+,J)$-free \then \, it
is $(\theta,J)$-stable.

\noindent
5) If $\cF \subseteq {}^\kappa \Ord$ is $(\theta^+_2,\theta_1,J)$-free
 \then \, $\cF$ is $\langle \theta_2,\theta_1,J\rangle$-free.

\noindent
6)  If $\cF \subseteq {}^\kappa \Ord$ is $\langle
   \theta_2,\theta_1,J\rangle$-free \then \, it is
   $(\theta^+_2,\theta_1,J)$-stable. 
\end{observation}

\begin{remark}
\label{b32}
We also have obvious monotonicity in $\cF$ and $\bar f$ and other
obvious implications.
\end{remark}

\begin{claim}
\label{b34}
1) Assume $\cF \subseteq {}^\kappa \Ord$ is semi-$(\theta,J)$-stable
   or just $J$ is $\theta_*$-complete and $\varepsilon \le \theta$.
\Then \, $\cF$ is semi-$(\theta^{+ \varepsilon +1},J)$-stable.

\noindent
2) Similarly without semi.
\end{claim}

\end{document}